\documentclass[10pt]{amsart}
\usepackage{amssymb}

\numberwithin{equation}{section}

\theoremstyle{plain}
\newtheorem{thm}[equation]{Theorem}
\newtheorem{cor}[equation]{Corollary}
\newtheorem{lem}[equation]{Lemma}
\newtheorem{prop}[equation]{Proposition}

\newtheorem{def-lem}[equation]{Lemma-Definition}

\theoremstyle{definition}
\newtheorem{defn}[equation]{Definition}

\theoremstyle{remark}
\newtheorem{rem}[equation]{Remark}
\newtheorem{ex}[equation]{Example}

\newcommand{\comment}[1]{}
\newcommand{\Ko}{K^\circ}
\newcommand{\Koo}{K^{\circ\circ}}

\newcommand{\Kt}{\widetilde K}

\newcommand{\Kalg}{K_{alg}}
\newcommand{\Kalgo}{K_{alg}^\circ}
\newcommand{\Kalgoo}{K_{alg}^{\circ\circ}}

\newcommand{\abs}[1]{\left\vert#1\right\vert}
\newcommand{\norm}[1]{\left\Vert#1\right\Vert_{\sup}}

\newcommand{\sq}{\,\square\,}
\newcommand{\ord}{\operatorname{ord}}

\renewcommand{\epsilon}{\varepsilon}
\renewcommand{\phi}{\varphi}
\renewcommand{\emptyset}{\varnothing}

\def\deg{{\rm deg}}
\def\ord{{\rm ord}}
\def\ac{{\overline{\rm ac}}}

\def\Var{{\rm Var}}

\def\Mod{{\rm Mod}}

\def\LPas{\cL_{\rm DP}}
\def\Lan{\cL_{S(A)}}
\def\Lans{\cL_{T(A)}}

\def\Tan{\mathcal{T}_{S(A)}}
\def\Tans{\mathcal{T}_{T(A)}}

\def\TPas{\mathcal{T}_{\rm DP}}

\def\Tanh{\mathcal{T}_{\mathrm{an}}^{\mathrm{h}}}
\def\Tanh0{\mathcal{T}_{{\rm an},0}^{\rm h}}

\def\Ord{{\rm Ord}}
\newcommand{\Res}{\mathrm{Res}}
\newcommand{\Val}{\mathrm{Val}}

\def\11{{\mathbf 1}}
\def\AA{{\mathbf A}}

\def\FF{{\mathbf F}}

\def\LL{{\mathbf L}}

\def\NN{{\mathbf N}}

\def\QQ{{\mathbf Q}}

\def\ZZ{{\mathbf Z}}

\def\cA{{\mathcal A}}
\def\cB{{\mathcal B}}

\def\cK{{\mathcal K}}
\def\cL{{\mathcal L}}
\def\cM{{\mathcal M}}

\def\cO{{\mathcal O}}

\def\cT{{\mathcal T}}
\def\cU{{\mathcal U}}

\begin{document}
\title[Analytic cell decomposition]{Analytic cell decomposition and
analytic motivic integration}
\author[Cluckers]{R.~Cluckers$^\mathrm{1}$}
\address{Katholieke Universit Leuven, Departement wiskunde,
Celestijnenlaan 200B, B-3001 Leu\-ven, Bel\-gium. Current address:
\'Ecole Normale Sup\'erieure, D\'epartement de
ma\-th\'e\-ma\-ti\-ques et applications, 45 rue d'Ulm, 75230 Paris
Cedex 05, France} \email{cluckers@ens.fr}
\urladdr{www.wis.kuleuven.ac.be/algebra/Raf/}

\thanks{$^\mathrm{1}$ The author
has been supported as a postdoctoral fellow by the Fund for
Scientific Research - Flanders (Belgium) (F.W.O.) during the
preparation of this paper.}
\thanks{$^\mathrm{2}$ The authors have been supported in part by NSF grants DMS-0070724 and DMS-0401175.}

\author[Lipshitz]{L.~Lipshitz$^\mathrm{2}$}
\address{Purdue University, West Lafayette IN 47907 USA}
\email{lipshitz@math.purdue.edu}

\author[Robinson]{Z.~Robinson$^\mathrm{2}$}
\address{East Carolina University, Greenville NC 27858 USA}
\email{robinsonz@mail.ecu.edu}

\subjclass[2000]{Primary 32P05, 32B05, 32B20, 03C10, 28B10;
Secondary 32C30, 03C98, 03C60, 28E99}

\begin{abstract}
The main results of this paper are a Cell Decomposition Theorem for
Henselian valued fields with analytic structure in an analytic
Denef-Pas language, and its application to analytic motivic
integrals and analytic integrals over $\FF_q((t))$ of big enough
characteristic. To accomplish this, we introduce a general framework
for Henselian valued fields $K$ with analytic structure, and we
investigate the structure of analytic functions in one variable,
defined on annuli over $K$. We also prove that, after
parameterization, definable analytic functions are given by terms.
The results in this paper pave the way for a theory of
\emph{analytic} motivic integration and \emph{analytic} motivic
constructible functions in the line of R.~Cluckers and F.~Loeser
[\emph{Fonctions constructible et int\'egration motivic I}, Comptes
rendus de l'Acad\'emie des Sciences, {\bf 339} (2004) 411 - 416].
\end{abstract}

\maketitle

\tableofcontents

\subsection{R\'esum\'e. D\'ecomposition cellulaire analytique et int\'egration analytique motivique. }
Dans cet article nous \'etablissons une d\'ecomposition cellulaire
pour des corps valu\'es Henseliens muni d'une structure analytique
induite par un langage de Denef-Pas analytique. En particulier, nous
appliquons cet \'enonc\'e \`a l' \'etude des int\'egrales
analytiques motiviques et des intégrales analytiques sur
$\FF_q((t))$ de charact\'eristic assez grand. Pour cela, il est
n\'ecessaire d'introduire une d\'efinition g\'en\'eral des corps
valu\'es Henselien $K$ avec structure analytique.  On examine alors
la structure de fonctions analytiques dans une variable d\'efinies
sur des annuli sur $K$ et l'on \'etablit que, dans ce contexte, les
fonctions d\'efinissables sont exactement donn\'es par des termes
apr\`es param\'etrisation. Plus g\'en\'eralement, les r\'esultats de
cet article pr\'eparent le chemin pour d\'efinir une th\'eorie
d'int\'egration analytique motivique et des fonctions analytiques
motiviques constructibles dans l'esprit de R.~Cluckers et F.~Loeser
[\emph{Fonctions constructible et int\'egration motivic I}, Comptes
rendus de l'Acad\'emie des Sciences, {\bf 339} (2004) 411 - 416].

\section{Introduction}

The main results of this paper are a Denef-Pas Cell Decomposition
Theorem for Henselian valued fields with analytic structure and its
application to analytic motivic integrals\footnote{Motivic here
stands for the idea of giving a geometric meaning to $p$-adic
integrals, uniform in $p$.} and analytic integrals over $\FF_q((t))$
of big enough characteristic. In section \ref{part:an} we introduce
a framework for Henselian valued fields $K$ with analytic structure,
both strictly convergent and separated, that generalizes \cite{bgr},
\cite{vdd1}, \cite{dhm}, \cite{LRSep}, \cite{unifrid}, and that
works in all characteristics. An analytic structure is induced by
rings of power series over a Noetherian ring $A$ that is complete
and separated with respect to the $I$-adic topology for some ideal
$I$ of $A$. This framework facilitates the use of standard
techniques from the theory of analytic rings in a more general model
theoretic setting.

\par
Another necessary ingredient for cell decomposition is an analysis
of analytic functions in one variable, defined on annuli over $K$.
This is carried out in Section~\ref{part:functions}.
Theorem~\ref{thm:torational} relates such functions piecewise to a
(strong) unit times a quotient of polynomials. This result is
extended to functions in one variable given by terms in
Theorem~\ref{lem:toholom}. The results of this section extend the
analysis of the ring of analytic functions on an affinoid subdomain
of $K$, carried out in \cite[Sections~2.1 and~2.2]{FP}, in the case
that $K$ is algebraically closed and complete, in three directions:
(i) $K$ not necessarily algebraically closed, (ii) the case of
quasi-affinoid subdomains of $K$ and (iii) the case that $K$ is not
complete but carries an analytic structure. This analysis will be
pursued further in a forthcoming paper.

\par
We also prove a fundamental structure result on definable analytic
functions, namely, that any definable function is given, after
parameterization using auxiliary sorts, by terms in a somewhat
bigger language, cf.~Theorem \ref{normal:an}. This structure
result is new also in the algebraic case, and is used in
\cite{CLoes} to prove a change of variables formula for motivic
integrals.

\par
The results in this paper pave the way for a theory of analytic
motivic integration and analytic motivic constructible functions in
the line of \cite{CLoes}, \cite{CLoes1}, \cite{CLoes2}. Another
approach to analytic motivic integration, based on entire models of
rigid varieties and the theory of N\'eron models instead of cell
decomposition, is developed by Sebag and Loeser in \cite{SB1}, and
by Sebag in \cite{SB2}. This alternative approach is pursued in
\cite{Seb2} for the study of generating power series and in
\cite{NS} for the study of the monodromy conjecture. In
\cite{CLoes}, apart from cell decomposition, also a dimension theory
is used; in the analytic case this can be developed along the lines
of work by \c{C}elikler \cite{Celikler}.

\subsection[]{} Cell decomposition is a technique of breaking
definable sets into finitely many definable pieces each of which
is particularly simple in a chosen coordinate direction. For
example, in the real case, Fubini's Theorem often reduces the
computation of an integral over a complicated set to an iterated
integral over the region between two graphs, on which the
integrand is of a simple form with respect to this coordinate,
cf.~the Preparation Theorem and its use for integration by Lion
and Rolin in \cite{LR}.

\par
In \cite{Cohen}, Cohen reproved Tarski's real quantifier elimination
using his real cell decomposition for definable sets. In the same
paper, he gave a cell decomposition for some Henselian fields,
e.g.~$p$-adic fields, extending results of Ax and Kochen. A cell
over a real field is a set given by conditions of the form
$f(x)<y<g(x)$ or $y=f(x)$, where $f,g$ are definable. That
quantifier elimination follows from cell decomposition is fairly
clear; the other implication a bit more complicated. A cell over a
Henselian field is specified by simple conditions on the order and
angular component of $y-c(x)$, where $c$ is definable (see below for
definitions). This reflects the idea that for many Henselian fields,
a statement about the field can be reduced to statements about the
value group and the residue field.

\par
Denef \cite{Denef2} refined Cohen's techniques to reprove
Macintyre's quantifier elimination for $p$-adic fields and to obtain
a $p$-adic integration technique which he used to prove the
rationality of certain $p$-adic Poincar\'e series \cite{Denef}.
 Pas \cite{pas}, \cite{pas2}, and Macintyre \cite{MacintRat}
extended this method to study uniform properties of $p$-adic
integrals.
 Denef and van den Dries \cite{DvdD} extended the
 Ax-Kochen-Cohen-Macintyre
$p$-adic quantifier elimination to the analytic category based on
strictly convergent power series. (See also \cite{vdd1}.)
 These ideas were extended to the algebraically closed analytic
category using separated power series by the second and third
authors, see \cite{mod}.
 The first author \cite{c}, using work of Haskell, Macpherson, and
van den Dries \cite{dhm}, obtained an analytic variant of the
$p$-adic cell decomposition and an application to $p$-adic
analytic integrals.

\par
 In this paper, we extend the ideas of quantifier elimination and
cell decomposition to a wider class of Henselian fields with
analytic structure, cf.~Theorems \ref{thm:qe:an} and \ref{thm:cd}.

\subsection{} Let us elaborate on the application to analytic
motivic integrals. Let $\cA$ be the class of all fields $\QQ_p$ for
all primes $p$ together with all their finite field extensions and
let $\cB$ be the class of all the fields $\FF_q((t))$ with $q$
running over all prime powers. For each fixed prime $p$ and integer
$n>0$ let $\cA_{p,n}$ be the subset of $\cA$ consisting of all
finite field extensions of $\QQ_p$ with degree of ramification fixed
by $\ord_p(p)=n$. For $K\in \cA\cup\cB$ write $K^\circ$ for the
valuation ring, $\widetilde K$ for the residue field, $\pi_K$ for a
uniformizer of $K^\circ$, and $q_K$ for $\sharp \widetilde K$.

\par
Denote by $\ZZ[[t]]\langle x_1,\ldots,x_n\rangle$ the ring of
strictly convergent power series over $\ZZ[[t]]$ (consisting of all
$\sum_{i\in\NN^n} a_i(t)x^i$ with $a_i(t)\in\ZZ[[t]]$ such that for
each $n\geq 0$ there exists $n'$ such that  $a_i(t)$ belongs to
$(t^n)$ for each $i$ with $i_1+\ldots+i_n>n'$).

\par
The purpose of these strictly convergent power series is to
provide analytic functions in a uniform way, as follows. To each $
f(x)=\sum_{i\in\NN^n} a_i(t)x^i$ in $\ZZ[[t]]\langle
x_1,\ldots,x_n\rangle$, $a_i(t)\in \ZZ[[t]]$, we associate for
each $K\in\cA\cup \cB$ the $K$-analytic function
   $$
   f_{K}:(K^\circ)^n\to K^\circ: x\mapsto \sum_{i\in\NN^n} a_i(\pi_K)
   x^i.
   $$
Fix $f\in \ZZ[[t]]\langle x_1,\ldots,x_n\rangle$. As $K$ varies in
$\cA\cup\cB$, one has a family of numbers
\begin{equation}\label{eq:ak}
a_K:=\int\limits_{(K^\circ)^n}|f_{K}(x)||dx|,
\end{equation}
and one would like to understand the dependence on $K$ in a
geometric way (see \cite{Kazhdan}
 for a context of this question).

\par
This is done using ${\rm Var}_\ZZ$, the collection of isomorphism
classes of algebraic varieties over $\ZZ$ (i.e. reduced separated
schemes of finite type over $\ZZ$), and, ${\rm Form}_\ZZ$, the
collection of equivalence classes of formulas in the language of
rings\footnote{Two formulas are equivalent in this language if they
have the same $R$-rational points for every ring $R$.} with
coefficients in $\ZZ$. For each finite field $k$, we consider the
ring morphisms
$$
{\rm Count}_{k}: \QQ[{\rm
Var}_\ZZ,\frac{1}{\AA^1_\ZZ},\{\frac{1}{1-\AA^{i}_\ZZ}\}_{i<0}
]\to\QQ
$$
which sends $Y\in{\rm Var}_\ZZ$ to $\sharp Y(k)$, the number of
$k$-rational points on $Y$, and,
$$
{\rm Count}_{k}: \QQ[{\rm
Form}_\ZZ,\frac{1}{\AA^1_\ZZ},\{\frac{1}{1-\AA^{i}_\ZZ}\}_{i<0}
]\to\QQ
$$
which sends $\varphi\in{\rm Form}_\ZZ$ to $\sharp \varphi(k)$, the
number of $k$-rational points on $\varphi$, and where we also write
$\AA^\ell_\ZZ$ for the isomorphism class of the formula
$x_1=x_1\wedge\ldots\wedge x_\ell=x_\ell$ (which has the set
$R^\ell$ as $R$-rational points for any ring $R$), $\ell\geq 0$.

\par
Using the work established in this paper, as well as results of
Denef and Loeser \cite{DL}, we establish Theorem \ref{thm:int:mets},
which is a generalization of the following.
\begin{thm}\label{thm:int}
\item[(i)] There exists a (non-unique) element
$$
X\in \QQ[{\rm
Var}_\ZZ,\frac{1}{\AA^1_\ZZ},\{\frac{1}{1-\AA^{i}_\ZZ}\}_{i<0} ]
   $$
and a number $N$ such that for each field $K\in\cA\cup\cB$ with
${\rm Char} \, \widetilde K>N$, one has
$$
a_K= {\rm Count}_{\widetilde K}(X).
$$
In particular, if ${\rm Char} \widetilde K>N$, then $a_K$ only
depends on $\widetilde K$.
 \item[(ii)]
For fixed prime $p$ and $n>0$, there exists a (non-unique) element
$$
X_{p,n}\in \QQ[{\rm
Form}_\ZZ,\frac{1}{\AA^1_\ZZ},\{\frac{1}{1-\AA^{i}_\ZZ}\}_{i<0} ]
   $$
such that for each field $K\in\cA_{p,n}$ one has
$$
a_K= {\rm Count}_{\widetilde K}(X_{p,n}).
$$
\end{thm}

\par
To prove Theorems \ref{thm:int} and \ref{thm:int:mets} we calculate
the $a_K$ by inductively integrating variable by variable, in a
uniform way, using analytic cell decomposition. By such
decomposition, one can partition the domain of integration uniformly
in $K\in\cA\cup\cB$, for big enough residue field characteristic,
and prepare the integrand on the pieces in such a way that the
integral with respect to a special variable becomes easy.

\par
There is possibly an alternative approach to prove Theorems
\ref{thm:int} and \ref{thm:int:mets} by using analytic embedded
resolutions of $f=0$ over (a ring of finite type over) $\ZZ[[t]]$,
if such a resolution exists. We do not pursue this approach.

\par
We comment on the non uniqueness of $X$  and $X_{p,n}$ in Theorem
\ref{thm:int}. By analogy to \cite{CLoes}, \cite{DL} and using the
results of this paper, one could associate unique objects, a motivic
integral, to the data used to define $a_K$. Such objects would live
in some quotient of $\QQ[{\rm
Form}_\ZZ,\frac{1}{\AA^1_\ZZ},\{\frac{1}{1-\AA^{i}_\ZZ}\}_{i<0} ]$,
and the morphisms ${\rm Count}_{k}$ could factor through this
quotient (at least for char $k$ big enough). To establish uniqueness
in some ring is beyond the scope of the present paper.

\par
More generally, we consider
\begin{equation}\label{eq:aks:1}
b_K(s):=\int\limits_{(K^\circ)^n}|f_{1K}(x)|^s|f_{2K}(x)||dx|,
\end{equation}
and similar integrals, when $K$ varies over $\cA\cup\cB$, where
$f_1,f_2$ are in $\ZZ[[t]]\langle x_1,\ldots,x_n\rangle$, and $s\geq
0$ is a real variable. In the generalization Theorem
\ref{thm:int:mets} of Theorem \ref{thm:int}, we prove the
rationality of $b_K(s)$ in $q_K^s$ for $K\in\cB$ of characteristic
big enough and with $q_K$ the number of elements in the residue
field of $K$. For $K\in \cA$ this rationality was proved by Denef
and van den Dries in \cite{DvdD}. This rationality is related to the
counting of points modulo $t^n$ in $\FF_q[[t]]/(t^n)$ for all $n$,
satisfying analytic equations, as is well known by the work of Igusa
and Denef.

\section{Analytic structures}\label{part:an}

Analytic structures, introduced in \cite{vdd1} (cf.~\cite{dhm} and
\cite{unifrid}), are a framework for the model theory of analytic
functions. This section contains an extensive elaboration of those
ideas.

\par
Model theory provides a convenient means to analyze algebraic
properties that depend on parameters, and analytic structures are
a way to extend model-theoretic techniques to the analytic
setting. In particular, a cell decomposition for a family of
functions of several variables is a partition of the domain into
finitely many simple sets on each of which the behavior of the
functions has a simple dependence on the value of the last
variable. By assigning the other variables fixed values in a
possibly non-standard field extension, the compactness theorem
reduces the problem of cell decomposition for polynomial functions
of several variables to that of obtaining a cell decomposition for
polynomial functions of just one variable, at the expense of
providing a uniform cell decomposition for all models of the
theory. By expressing a power series as the product of a unit and
a polynomial in the last variable, Weierstrass Preparation is used
to reduce analytic questions to algebraic ones. An analytic
structure provides a convenient framework for dealing with the
parameters that arise in applying Weierstrass Preparation.

\par
To make use of the Weierstrass data, the definition of a
model-theoretic structure must be extended so that compositional and
algebraic identities in the power series ring are preserved when the
power series are interpreted as functions on the underlying field.
In the case of polynomial rings, the interpretation of addition and
multiplication in a model of the theory of rings already provides a
natural homomorphism from the polynomial ring into the ring of
functions on the underlying structure. Furthermore, if the
underlying field is complete, the valued field structure also
already provides a natural homomorphism from the ring of convergent
power series into the ring of functions (that preserves not only
algebraic, but compositional identities as well). But the fields
over which we work may not be complete since we must work uniformly
in all models of a given theory. Thus, to apply the Weierstrass
techniques, our models must come equipped with a distinguished
homomorphism from the ring of power series to the ring of functions.
This is essentially the definition of analytic structure in
Definition~\ref{analytic_structure} and in \cite{vdd1}. (Note that,
rather than a distinguished homomorphism, one could employ instead a
first-order axiom scheme in which each power series identity is
coded into an axiom, but that obscures the difference between the
algebraic and analytic situations, where topological completeness,
in some form, comes into play.)

\par
As in \cite{vdd1}, \cite{dhm} and \cite{unifrid}, in using
Weierstrass techniques, one often introduces new parameters for
certain ratios. Without a natural means of adjoining elements of a
(possibly non-standard) model to the given coefficient ring of the
power series ring, one is prevented from specializing the
parameters, which complicates some computations. However, since the
proof of the Weierstrass Division Theorem relies on completeness in
the coefficient ring, adjoining elements of an arbitrary model to
the coefficient ring is problematic. The methods of \cite{LRSep}
were developed to analyze the commutative algebra of rings of
separated power series, which are filtered unions of complete rings.
Those ideas are applied in this section to show how to extend the
coefficient ring (Theorem~\ref{sigma_completion} and
Definition~\ref{power_series_with_K_parms}) and ground field
(Theorem~\ref{algebraic_extension}) of a given analytic structure,
which is how the present treatment of analytic structures differs
from the previous ones. (Indeed, with minor modifications to the
proofs, much of the theory of \cite{LRSep} applies to the rings
$S_{m,n}(\sigma,K)$ introduced in
Definition~\ref{power_series_with_K_parms}, and, although we prefer
to give a self-contained treatment in this paper, would simplify the
proofs of the results of Section~\ref{part:functions}.)

\par
Finally, let $K$ be a separated analytic $A$-structure as in
Definition~\ref{analytic_structure}, so the power series in a ring
$S_{m,n}(A)$ are interpreted as analytic functions on $K$ in such a
way as to preserve the algebraic and compositional identities of
$S_{m,n}(A)$  and an extended power series ring $S_{m,n}(\sigma,K)$
is obtained from $S_{m,n}(A)$, as in
Definition~\ref{power_series_with_K_parms}, by adjoining
coefficients from the field $K$. It is important to note that,
although the extended power series rings $S_{m,n}(\sigma,K)$ are
much larger than the rings $S_{m,n}(A)$, the structure $K$ has
essentially the same first-order diagram in the extended language.
Thus, although it is easier to work with the extended power series
rings $S_{m,n}(\sigma,K)$, they have the same model theory as the
smaller rings $S_{m,n}(A)$, which, in fact, is the point of
introducing the extension.

\begin{defn}
\label{strictly_convergent} Let $E$ be a Noetherian ring that is
complete and separated for the $I$-adic topology, where $I$ is a
fixed ideal of $E$. Let $(\xi_1,\dots,\xi_m)$ be variables. The
ring of {\em strictly convergent power series in $\xi$ over $E$}
(cf.~\cite{bgr}, Section~1.4) is
$$
T_m(E)=E\langle\xi\rangle:=\{ \sum_{\nu\in\mathbb
N^m}a_\nu\xi^\nu\colon \lim_{|\nu|\to\infty}a_\nu=0\}.
$$
Let $(\rho_1,\dots,\rho_n)$ be variables. The ring
$$
S_{m,n}(E):=E\langle\xi\rangle[\![\rho]\!]
$$
is a ring of {\em separated power series over $E$}
(cf.~\cite{LRSep}, Section~2).
\end{defn}

\begin{rem}
\item[(i)]
If the formal power series ring $E[\![\rho]\!]$ is given the
ideal-adic topology for the ideal generated by $I$ and $(\rho)$,
then $S_{m,n}(E)$ is isomorphic to
$E[\![\rho]\!]\langle\xi\rangle$. Note that $S_{m,0}=T_m$ and
$T_{m+n}(E)$ is contained in $S_{m,n}(E)$.

\item[(ii)]
\label{completion_of_polynomials} Observe that
$E\langle\xi\rangle$ is the completion of the polynomial ring
$E[\xi]$ in the $I\cdot E[\xi]$-adic topology and
$E\langle\xi\rangle[\![\rho]\!]$ is the completion of the
polynomial ring $E[\xi,\rho]$ in the $J$-adic topology, where $J$
is the ideal of $E[\xi,\rho]$ generated by $\rho$ and the elements
of $I$.

\item[(iii)] The example of $A=\ZZ[[t]]$ and $I=(t)$ is the one used
in the introduction to put a strictly convergent analytic structure
on the $p$-adic fields and on the fields $\FF_q((t))$.
\end{rem}

The Weierstrass Division Theorem (cf.~\cite{LRSep}, Theorems~2.3.2
and~2.3.8) provides a key to the basic structure of the power
series rings $S_{m,n}(E)$.

\begin{defn}
\label{regular} Let $f\in S_{m,n}(E)$. The power series $f$ is
{\em regular in $\xi_m$ of degree $d$} if, and only if, $f$ is
congruent, modulo the ideal $I+(\rho)$, to a monic polynomial in
$\xi_m$ of degree $d$, and $f$ is {\em regular in $\rho_n$ of
degree $d$} if, and only if, $f$ is congruent, modulo the ideal
$I+(\rho_1,\dots,\rho_{n-1})$ to $\rho_n^d\cdot g(\xi,\rho)$ for
some unit $g$ of $S_{m,n}(E)$.
\end{defn}

\begin{prop}[Weierstrass Division]
\label{weierstrass_division} Let $f,g\in S_{m,n}(E)$.

{\em (i)}  Suppose that $f$ is regular in $\xi_m$ of degree $d$.
Then there exist uniquely determined elements $q\in S_{m,n}(E)$
and $r\in S_{m-1,n}(E)[\xi_m]$ of degree at most $d-1$ such that
$g=qf+r$. If $g\in J\cdot S_{m,n}$ for some ideal $J$ of
$S_{m-1,n}$, then $q,r\in J\cdot S_{m,n}$.

{\em (ii)} Suppose that $f$ is regular in $\rho_n$ of degree $d$.
Then there exist uniquely determined elements $q\in S_{m,n}(E)$
and $r\in S_{m,n-1}(E)[\rho_n]$ of degree at most $d-1$ such that
$g=qf+r$. If $g\in J\cdot S_{m,n}$ for some ideal $J$ of
$S_{m,n-1}$, then $q,r\in J\cdot S_{m,n}$.
\end{prop}

\begin{rem}
\label{strictly_convergent_weierstrass}
By taking $n=0$ in Proposition~\ref{weierstrass_division}~(i), one
obtains a Weierstrass Division Theorem for the $T_m(E)$.
\end{rem}

Dividing $\xi_m^d$ (respectively, $\rho_n^d$) by an element
$f\in S_{m,n}$ regular in $\xi_m$, respectively, $\rho_n$, of degree
$d$, as in \cite{LRSep}, Corollary~2.3.3, we obtain the following
corollary.

\begin{cor}[Weierstrass Preparation]
\label{weierstrass_preparation} Let $f\in S_{m,n}(E)$.

{\em (i)} If $f$ is regular in $\xi_m$ of degree $d$, then there
exist: a unique unit $u$ of $S_{m,n}$ and a unique monic
polynomial $P\in S_{m-1,n}[\xi_m]$ of degree $d$ such that
$f=u\cdot P$.

{\em (ii)} If $f$ is regular in $\rho_n$ of degree $d$, then there
exist: a unique unit $u$ of $S_{m,n}$ and a unique monic
polynomial $P\in S_{m,n-1}[\rho_n]$ of degree $d$ such that
$f=u\cdot P$; in addition, $P$ is regular in $\rho_n$ of degree
$d$.
\end{cor}

Let the ring $E$ and the ideal $I$ be as in
Definition~\ref{strictly_convergent}. If $K$ is a field containing
$E$ that is complete in a rank~1 valuation and $I$ is contained in
the maximal ideal $\Koo$ of the valuation ring $\Ko$, then
$T_m(E)$ (respectively, $S_{m,n}(E)$) may be interpreted as a ring
of analytic functions on the polydisc $(\Ko)^m$ (respectively,
$(\Ko)^m\times(\Koo)^n$), exactly as in \cite{LRSep}. The
following definition permits an extension to more general valued
fields $K$, for example, to ultraproducts of complete fields. In
this more general setting, analytic properties usually derived
employing the completeness of the domain can often be derived
instead from Weierstrass Division (which relies on completeness in
the coefficient ring).

\begin{defn}[cf.~\cite{vdd1} and \cite{unifrid}]
\label{analytic_structure} Let $A$ be a Noetherian ring that is
complete and separated with respect to the $I$-adic topology for a
fixed ideal $I$ of $A$. Let $(K,\ord,\Gamma)$ be a valued field. A
{\em separated analytic $A$-structure on $K$} is a collection of
homomorphisms $\sigma_{m,n}$ from $S_{m,n}(A)$ into the ring of
$\Ko$-valued functions on $(\Ko)^m\times(\Koo)^n$ for each $m,n\in
\mathbb N$ such that:
\begin{enumerate}
\item[(i)] $(0)\not= I\subset \sigma_0^{-1}(\Koo)$,

\item[(ii)] $\sigma_{m,n}(\xi_i)=$ the $i$-th coordinate function
on $(\Ko)^m\times(\Koo)^n$, $i=1,\dots,m$, and
$\sigma_{m,n}(\rho_j)=$ the $(m+j)$-th coordinate function on
$(\Ko)^m\times(\Koo)^n$, $j=1,\dots,n$,

\item[(iii)] $\sigma_{m,n+1}$ extends $\sigma_{m,n}$, where we
identify in the obvious way functions on $(\Ko)^m\times(\Koo)^n$
with functions on $(\Ko)^m\times(\Koo)^{n+1}$ that do not depend
on the last coordinate, and $\sigma_{m+1,n}$ extends
$\sigma_{m,n}$ similarly.
\end{enumerate}
A collection of homomorphisms $\sigma_m$ from $T_m(A)=S_{m,0}(A)$
into the ring of $\Ko$-valued functions on $(\Ko)^m$ is called a
{\em strictly convergent analytic $A$-structure on $K$} if the
homomorphisms $\sigma_{m,0}:=\sigma_m$ satisfy the above three
conditions (with $n=0$).
\end{defn}

Here are some typical examples of valued fields with strictly
convergent analytic $A$-structure. Take $A:=\mathbb Z[\![t]\!]$,
where $t$ is one variable, equipped with the $(t)$-adic topology.
Then $(\mathbb C(\!(t)\!),\ord_t,\mathbb Z)$ carries a unique
analytic $A$-structure determined by $\sigma_0(t)=t$. For each prime
$p\in\mathbb N$, the valued field of $p$-adic numbers $(\mathbb
Q_p,\ord_p,\mathbb Z)$ carries a unique analytic $A$-structure
determined by $\sigma_0(t)=p$. Similarly, $\sigma_0(t)=p$ determines
a unique separated analytic structure on the non-discretely valued
field $\mathbb C_p$, the completion of the algebraic closure of
$\mathbb Q_p$, which yields a larger family of analytic functions
than the corresponding strictly convergent analytic $A$-structure.
The fields $\FF_p((t))$ carry unique analytic $A$-structures
determined by $\sigma_0(t)=t$. The latter (standard) analytic
$A$-structures induce an analytic $A$-structure on any non-principal
ultraproduct of the $p$-adic fields $\mathbb Q_p$, or $\mathbb C_p$,
or $\mathbb F_p((t))$. Note that such fields carry analytic
$A$-structure even though they are not complete.

\par
By definition, analytic $A$-structures preserve the ring operations
on power series, thus they preserve the Weierstrass Division data.
It follows that analytic $A$-structures also preserve the operation
of composition.

\begin{prop}
\label{composition} Analytic $A$-structures preserve composition.
More precisely, if $f\in S_{m,n}(A)$, $\alpha_1,\dots,\alpha_m\in
S_{M,N}(A)$, $\beta_1,\dots,\beta_n\in IS_{M,N}(A) +(\rho)$, where
$S_{M,N}(A)$ contains power series in the variables $(\xi,\rho)$
and $I$ is the fixed ideal of $A$, then $g:=f(\alpha,\beta)$ is in
$S_{M,N}(A)$ and
$\sigma(g)=(\sigma(f))(\sigma(\alpha),\sigma(\beta))$.
\end{prop}

\begin{proof} By the Weierstrass Division
Theorem, there are elements $q_i\in S_{m+M,n+N}(A)$ such that
$$
f(\eta,\lambda)=g(\xi,\rho)+\sum_{i=1}^m(\eta_i-\alpha_i(\xi,\rho))\cdot
q_i+ \sum_{j=1}^n(\lambda_j-\beta_j(\xi,\rho))\cdot q_{m+j}.
$$
Let $(x,y)\in(\Ko)^M\times(\Koo)^N$ and put
$a_i:=\sigma(\alpha_i)(x,y)$ and $b_j:=\sigma(\beta_j)(x,y)$ for
$i=1,\dots,m$, $j=1,\dots,n$. Clearly, $b\in(\Koo)^n$. By plugging
$(a,x,b,y)$ into the above equation, the proposition follows.
\end{proof}

Next we show that the image of a power series is the zero function
if, and only if, the image of each of its coefficients is zero.
This employs parameterized Weierstrass Division which relies on
the strong Noetherian property of Lemma~\ref{strong_noetherian}.

\begin{lem}[\cite{LRSep}, Lemma~3.1.6]
\label{strong_noetherian} Let $F\in S_{m+M,n+N}(E)$ and write
$$
F=\sum f_{\mu,\nu}(\xi,\rho)\eta^\mu\lambda^\nu
$$
for some
$f_{\mu,\nu}\in S_{m,n}(E)$. Then there are: $d\in\mathbb N$ and
units $G_{\mu,\nu}$ of $S_{m+M,n+N}(E)$ such that
$$
F=\sum_{|\mu|+|\nu|\le d}f_{\mu,\nu}(\xi,\rho)\eta^\mu\lambda^\nu
G_{\mu,\nu}(\xi,\eta,\rho,\lambda).
$$
\end{lem}

\begin{prop}\label{kernel}
The image of a power series is the zero function if, and only if,
the image of each of its coefficients is zero. More precisely,

{\em (i)} Let $\sigma$ be a separated analytic $A$-structure on
the valued field $K$. Then $\ker\sigma_{m,n}=\ker\sigma_0\cdot
S_{m,n}(A)$.

{\em (ii)} Let $\sigma$ be a strictly convergent analytic
$A$-structure on the valued field $K$. Then
$\ker\sigma_m=\ker\sigma_0\cdot T_m(A)$.
\end{prop}

\begin{proof}
(ii) The ring $\sigma_0(A)$ is a Noetherian ring that is complete
and separated in the $\sigma_0(I)$-adic topology. The map $\sigma$
induces a homomorphism $\pi_m\colon T_m(A)\to T_m(\sigma_0(A))$, and
$\ker\pi_m=\ker\sigma_0\cdot T_m(A)$. Thus, the homomorphism
$\sigma_m$ factors through $T_m(\sigma_0(A))$, yielding a strictly
convergent $\sigma_0(A)$-analytic structure $\bar\sigma$ on $\Ko$.
Hence, there is no loss in generality to assume that
$\ker\sigma_0=(0)$. Let $f\in T_m(A)\setminus\{0\}$; we must show
that $\sigma_m(f)$ is not the zero function.
\par
Observe that if $f\in T_{m-1}[\xi_m]$ is monic in $\xi_m$ then
$\sigma_m(f)$ is not the zero function. Indeed, write
$f=\xi_m^d+\sum_{i=0}^{d-1}\xi_m^ia_i(\xi')$, where
$\xi'=(\xi_1,\dots,\xi_{m-1})$. Let $x\in(\Ko)^{m-1}$; then
$\sigma_m(f)(x,\xi_m)\in\Ko[\xi_m]$ is monic of degree $d$. By
Definition~\ref{analytic_structure}~(i), $K$ is a non-trivially
valued hence infinite field; thus $\sigma_m(f)$ is not the zero
function.

\par
Now let $f=\sum f_\mu\xi^\mu$ be any non-zero element of $T_m$. By
Lemma~\ref{strong_noetherian}, there are $d\in\mathbb N$ and units
$g_\mu$ of $T_m$ such that
$$
f=\sum_{|\mu|\le d}f_\mu\xi^\mu g_\mu.
$$
Let $\nu$ be the lexicographically largest index such that
$$
\ord\sigma(f_\nu)=\min_{|\mu|\le d}\ord\sigma(f_\mu).
$$
Let $\eta_\mu$ be new variables and put
$$
F:=\xi^\nu+\mathop{\sum_{|\mu|\le d}}_{\mu\neq\nu}\eta_\mu\xi^\mu
g_\nu^{-1}g_\mu.
$$
Since the above sum is finite,
$$
\sigma(f)=\sigma(f_\nu g_\nu)\sigma(F)(\xi,y_\mu),
$$
where the $y_\mu:=\frac{\sigma(f_\mu)}{\sigma(f_\nu)}\in\Ko$.
Since $g_\nu$ is a unit and $\sigma(f_\nu)\neq0$, to show that
$\sigma(f)$ is not the zero function, it suffices to show that
$\sigma(F)(\xi,y_\mu)$ is not the zero function.

\par
By the choice of $\nu$, there is a polynomial change of variables
$\phi$, involving only the $\xi$, such that $F\circ\phi$ is regular
in $\xi_m$ of some degree $d$. By Proposition~\ref{composition}, it
is enough to show that $\sigma(F\circ\phi)$ is not the zero
function, which follows from
Corollary~\ref{weierstrass_preparation}~(i).

\par
(i) As above, we may assume that $\ker\sigma_0=(0)$. Let
$f(\xi,\rho)$ be a nonzero element of $S_{m,n}$, let $a$ be a
nonzero element of the ideal $I$ of $A$ and put
$g:=f(\xi,a\cdot\rho)$. By Proposition~\ref{composition}, to show
that $\sigma(f)$ is not the zero function, it suffices to show
that $\sigma(g)$ is not the zero function. Since
$\ker\sigma_0=(0)$ and $\Ko$ is an integral domain, so is $A$.
Since $a$ is a nonzero element of $I$, $g$ is a nonzero element of
$T_{m+n}$. It is then a consequence of part~(ii) that $\sigma(g)$
is not the zero function.
\end{proof}

Next, we discuss how extend the coefficient ring of a given
analytic structure.

\begin{defn}
\label{compatible} Let $A$ and $E$ be Noetherian rings that are
complete and separated for the $I$-adic, respectively, $J$-adic,
topologies, where $I$ and $J$ are fixed ideals of $A$,
respectively, of $E$. Let $K$ be a valued field with analytic
$A$-structure $\{\sigma_{m,n}\}$ and analytic $E$-structure
$\{\tau_{m,n}\}$. Suppose $E$ is an $A$-algebra via the
homomorphism $\phi\colon A\to E$, and that $I\subset\phi^{-1}(J)$.
Note that $\phi$ extends coefficient-wise to a homomorphism
$\phi:S_{m,n}(A)\to S_{m,n}(E)$. The analytic structures  $\sigma$
and $\tau$ are called {\em compatible} if, for all $m$ and $n$,
$\sigma_{m,n}=\tau_{m,n}\circ\phi$.
\end{defn}

\par
It can be particularly useful to extend the coefficient ring of an
analytic $A$-structure by adjoining finitely many parameters from
the domain $K$. The coefficient rings of analytic structures are
complete, and Lemma~\ref{sigma_completion_existence} permits us to
define the appropriate completion of a finitely generated
$A$-subalgebra of $\Ko$.

\begin{def-lem}
\label{sigma_completion_existence} {\em (i)} Let $K$ be a valued
field with separated analytic $A$-structure $\{\sigma_{m,n}\}$ and
let $E$ be a finitely generated $\sigma_0(A)$-subalgebra of $\Ko$,
say, generated by $a_1,\dots,a_m$. Then $E$ is Noetherian. Let
$b_1,\dots,b_n$ generate the ideal $E\cap\Koo$. The subset
$E^\sigma$ of $K$
$$ E^\sigma:=\{\sigma(f)(a,b):
f\in S_{m,n}(A)\}
$$
is independent of the choices of $a$ and $b$. Moreover, $
E^\sigma$ is a Noetherian ring that is complete and separated with
respect to the $J$-adic topology, where $J=(E^\sigma\cap\Koo)$.
Moreover, $J$ is generated by $b_j$, $j=1,\dots,n$.

{\em (ii)} Let $K$ be a valued field with strictly convergent
analytic $A$-structure $\{\sigma_m\}$ and $E$ be a finitely
generated $\sigma_0(A)$-subalgebra of $\Ko$, say, generated by
$a_1,\dots,a_m$. Then $E$ is Noetherian. The subset $E^\sigma$ of
$K$
$$
E^\sigma:=\{\sigma(f)(a): f\in T_m(A)\}
$$
is independent of the choice of $a$. Moreover, $ E^\sigma$ is a
Noetherian ring that is complete and separated with respect to the
$J$-adic topology, where $J=\sigma_0(I)\cdot E^\sigma$.
\end{def-lem}

\begin{proof}
(i) Let $E$ be generated by some tuple $a'$ and $E\cap\Koo$ by
$b'$. For some polynomials $p_i,q_{j,\ell}\in A[\xi]$,
$$
a_i'=\sigma(p_i)(a),\ i=1,\dots,m' \quad\mathrm{and}\quad
b_j'=\sum_{\ell=1}^n\sigma(q_{j,\ell})(a)b_\ell,\ j=1,\dots,n'.
$$
That $E^\sigma$ is independent of the choice of $a$ now follows from
Proposition~\ref{composition}.

\par
To prove the remainder of part~(i), observe that the ideal $J$ of $
E^\sigma$ is generated by $b_1,\dots,b_n$. Indeed, let $f\in
S_{m,n}$ and write
$$
f=f_0(\xi)+\sum_{\ell=1}^r e_\ell g_\ell(\xi) +
\sum_{j=1}^n\rho_jh_j(\xi,\rho),
$$
where $f_0\in A[\xi]$ is a polynomial, the $e_\ell$ generate the
ideal $I$ of $A$, $g_\ell\in S_{m,0}$ and $h_j\in S_{m,n}$. Note
that $\ord\,\sigma_0(e_\ell),\ord\,\sigma(\rho_j)(b)=\ord\,b_j>0$ and
that $\sigma_0(e_\ell)$ and the $\sigma(\rho_j)(b)=b_j$ belong to
the ideal generated by $b$. Thus, $\ord\,\sigma(f)(a,b)>0$ implies
that $\ord\,\sigma(f_0)(a)>0$. Since $f_0$ is a polynomial,
$\sigma(f_0)(a)\in E$, and it follows that $\sigma(f_0)(a)$ must
also belong to the ideal generated by $b$.

\par
Now consider the $A$-algebra homomorphism
$$
\epsilon_{a,b}\colon S_{m,n}(A)\to E^\sigma\colon
f\mapsto\sigma(f)(a,b).
$$
Since $\epsilon_{a,b}$ is clearly surjective and $S_{m,n}(A)$ is
Noetherian, $ E^\sigma$ is Noetherian. By the above observation,
the non-trivial ideal $J$ is generated by the images of the
$\rho_j$ under $\epsilon_{a,b}$. Since $S_{m,n}$ is complete in
the $(\rho)$-adic topology, it follows from the Artin-Rees Theorem
that the finitely generated $S_{m,n}$-module $ E^\sigma$ is
complete and separated in the $J$-adic topology, as desired.

\par
(ii) The proof is similar to part~(i).
\end{proof}

Theorem~\ref{sigma_completion}, below, gives a basic example of
extending the coefficient ring of an $A$-analytic structure to
obtain a compatible analytic structure.

\begin{thm}
\label{sigma_completion} {\em (i)} Let $K$ be a valued field with
separated analytic $A$-structure $\{\sigma_{m,n}\}$. Let
$E\subset\Ko$ be a finitely generated $A$-subalgebra of $\Ko$ and
let $ E^\sigma$ be as in
Definition~\ref{sigma_completion_existence}~(i). Then $\sigma$
induces a unique analytic $ E^\sigma$-structure $\tau$ on $\Ko$
such that $\sigma$ and $\tau$ are compatible. Moreover, each
$\tau_{m,n}$ is injective.

{\em (ii)} The analogous statement holds for $K$ a valued field with
strictly convergent analytic $A$-structure $\{\sigma_m\}$.

\end{thm}

\begin{proof}
(i) Let $f\in S_{M,N}( E^\sigma)$. By
Lemma~\ref{sigma_completion_existence}, $J$ is generated by the
$\sigma(\rho_j)(b)$, so there is some $F\in S_{m+M,n+N}(A)$,
$F=\sum f_{\mu,\nu}(\xi,\rho)\eta^\mu\lambda^\nu$, such that
$$
f=\sum\sigma(f_{\mu,\nu})(a,b)\eta^\mu\lambda^\nu.
$$
Once the required homomorphisms $\tau$ are shown to exist, it
follows by the Weierstrass Division Theorem as in the proof of
Proposition~\ref{composition}, that
\begin{equation}
\label{assignment}
\tau_{m,n}(f)(\eta,\lambda)=\sigma_{m+M,n+N}(F)(a,\eta,b,\lambda);
\end{equation}
i.e., that $\tau_{m,n}$ is uniquely determined by the conditions
of Definition~\ref{compatible}.

\par
It remains to show that $\tau_{m,n}$ is well-defined by the
assignment of equation~\ref{assignment}. For that, it suffices to
show for any $G\in S_{m+M,n+N}(A)$, $G=\sum
g_{\mu,\nu}(\xi,\rho)\eta^\mu\lambda^\nu$, that if $\sum
g_{\mu,\nu}(a,b)\eta^\mu\lambda^\nu$ is the zero power series of
$S_{M,N}( E^\sigma)$, then $\sigma_{m+M,n+N}(G)(a,\eta,b,\lambda)$
is the zero function. By Lemma~\ref{strong_noetherian}, there are:
$d\in\mathbb N$ and power series $H_{\mu,\nu}\in S_{m+M,n+N}(A)$
such that
$$
G=\sum_{|(\mu,\nu)|\le d}g_{\mu,\nu}H_{\mu,\nu}.
$$
Then
$$
\sigma(G)(a,\eta,b,\lambda)= \sum_{|(\mu,\nu)|\le
d}\sigma(g_{\mu,\nu})(a,b)\sigma(H_{\mu,\nu})(a,\eta,b,\lambda)=0,
$$
as desired. Since $ E^\sigma$ is a subring of $\Ko$, the
injectivity of $\tau$ is a consequence of
Proposition~\ref{kernel}. This proves part~(i).

\par
(ii) The proof of part~(ii) is similar.
\end{proof}

For our purposes, it is useful to work with the ring of all
separated (or strictly convergent) power series with parameters
from $K$.

\begin{defn}
\label{power_series_with_K_parms} (i) Let $K$ be a valued field
with separated analytic $A$-structure $\{\sigma_{m,n}\}$. Let
$\mathcal F(\sigma,K)$ be the collection of all finitely generated
$A$-subalgebras $E\subset\Ko$. Then $\mathcal F(\sigma,K)$ and $\{
E^\sigma\}_{E\in\mathcal F(\sigma,K)}$ form direct systems of
$A$-algebras in a natural way, where $E^\sigma$ is as in
Definition~\ref{sigma_completion_existence}~(i). Put
$$
S_{m,n}^\circ(\sigma, K):=
\lim\limits_{\overrightarrow{E\in\mathcal F(\sigma,K)}}
E^\sigma\langle\xi\rangle[\![\rho]\!],
$$
which is a $\Ko$-algebra. The {\em rings of separated power series
with parameters from $K$} are then defined to be
$$
S_{m,n}(\sigma,K):= K\otimes_{\Ko}S_{m,n}^\circ(\sigma,K).
$$
(ii) Let $K$ be a valued field with strictly convergent analytic
$A$-structure $\{\sigma_m\}$. The rings $T_m(\sigma,K)$ of {\em
strictly convergent power series with parameters from $K$} are
defined similarly using
Lemma~\ref{sigma_completion_existence}~(ii).
\end{defn}

\begin{rem}
\label{weierstrass_remark} The rings $S_{m,n}^\circ(\sigma,K)$
(respectively, $T_m(\sigma,K)$) inherit Weierstrass Division,
Theorem~\ref{weierstrass_division}, and Weierstrass Preparation,
Corollary~\ref{weierstrass_preparation}, since they are direct
unions of the rings $S_{m,n}( E^\sigma)$ (respectively, $T_m(
E^\sigma)$) to which those results apply.
\end{rem}

Just as it can be useful to extend the coefficient ring of an
analytic structure, it is also useful to be able to extend the
domain of an analytic structure. This requires the following
proposition, which is proved exactly as \cite{cs}, Lemma~3.3.

\begin{prop}
\label{henselian} (i) Let $K$ be a valued field with separated
analytic $A$-structure; then $\Ko$ is a Henselian valuation ring.

(ii) Let $K$ be a valued field with strictly convergent analytic
$A$-structure such that $\ord(\Koo)$ has a minimal element
$\gamma$, and $\gamma=\min\ord(\sigma_0(I))$. Then $\Ko$ is a
Henselian valuation ring.
\end{prop}

The following theorem permits us to work over any finite algebraic
extension, or over the algebraic closure, of the domain of an
analytic $A$-structure.

\begin{thm}
\label{algebraic_extension} (i) Let $K$ be a valued field with
separated analytic $A$-structure $\sigma$. Then there is a unique
extension of $\sigma$ to a separated analytic $A$-structure $\tau$
on $\Kalg$, the algebraic closure of $K$.

(ii) Let $K$ be a valued field with strictly convergent analytic
$A$-structure such that $\ord(\Koo)$ has a minimal element
$\gamma$, and $\gamma=\min\ord(\sigma_0(I))$. Then there is a
unique extension of $\sigma$ to a strictly convergent analytic
$A$-structure $\tau$ on $\Kalg$.

(iii) Let $K$ be as in part~(ii); then there is a unique extension
of $\sigma$ to a separated analytic $A$-structure $\tau$ on $\Kalg$.
\end{thm}

\begin{proof}
Let $\alpha\in\Kalgo$ and let $P(t)=t^d+a_1t^{d-1}+\cdots+a_d$ be
the minimal polynomal for $\alpha$ over $K$. Since by
Proposition~\ref{henselian}~(i), $\Ko$ is Henselian, the
coefficients $a_i$ lie in $\Ko$; moreover, if $\alpha\in\Kalgoo$,
then the $a_i$ lie in $\Koo$. Now use Weierstrass division.
\end{proof}

\begin{rem}
\label{tensor}
Let $K$ be a valued field with analytic $A$-structure
that satisfies the conditions of either
Theorem~\ref{algebraic_extension}~(i) or~(ii), and let $L$ be an
extension of $K$ contained in $\Kalg$. Then the arguments of
Theorem~\ref{algebraic_extension} show that
$$
S_{m,n}^\circ(\tau,L)=
L^\circ\otimes_{\Ko}S_{m,n}^\circ(\sigma,K)\quad \mathrm{and}
\quad S_{m,n}(\tau,L)= L\otimes_K S_{m,n}(\sigma,K).
$$
\end{rem}

Since the base change is faithfully flat, Remark~\ref{tensor} yields the following corollary.

\begin{cor}
\label{flatness}
Let $K$ and $L$ be as in Remark~\ref{tensor}; then:

(i) $S_{m,n}^\circ(\tau,L)$ (respectively, $S_{m,n}(\tau,L)$) is
faithfully flat over $S_{m,n}^\circ(\sigma,K)$ (respectively, over
$S_{m,n}(\sigma,K)$), and

(ii) if $L$ is finite over $K$, then $S_{m,n}(\tau,L)$ is finite
over $S_{m,n}(\sigma,K)$.

Similar statements hold for $T_m$.
\end{cor}

\section{Rational analytic functions in one variable}\label{part:functions}

In this section, we develop the basis of a theory of analytic
functions on a $K$-annulus (an irreducible R-domain in  $K^\circ$),
when $K$ carries a separated $A$-analytic structure, as it is needed
for the proof of the cell decomposition of this paper. In
particular, we show that given an analytic function $f$ on a
$K$-annulus, there is a partition of the annulus into finitely many
annuli $\cU$ such that $f|_{\cU}$ is a rational function times a
(very) strong unit (see Theorem~\ref{thm:torational}). All the same
results hold (with the same proofs) in the ``standard'' case where
$K$ is a complete non-Archimedean valued field and
$S_{m,n}(\sigma,K)$ is replaced by $S_{m,n}(E,K)$. Hence, the
results in this section also extend the affinoid results of
\cite{FP}, Sections~2.1 and~2.2, to the case that $K$ is not
algebraically closed and to the quasi-affinoid case (i.e., allowing
strict as well as weak inequalities).

\par
The results of this section require $K$ to carry a separated
$A$-analytic structure. Note, however, by
Theorem~\ref{algebraic_extension}~(iii), in the setting of this
paper, a strictly convergent $A$-analytic structure on $K$ can be
extended uniquely to a separated $A$-analytic structure on $\Kalg$.

\par
A subsequent paper will give a complete treatment of the analytic
geometry of the one-dimensional unit ball over $\Kalg$, when $K$
carries either a strictly convergent or separated analytic
structure. This will include the analogue of the classical
Mittag-Leffler Theorem (cf. \cite{FP}, Theorems~2.2.6 and~2.2.9)
over coefficient fields $K$ that may be neither complete nor
algebraically closed, both in the affinoid and quasi-affinoid
setting. This will allow the exploration of more cell
decompositions.

\begin{defn}
\label{Rational definitionA}
Let $K$ be a Henselian
valued field (with separated $A$-analytic structure).

(a) A {\em $K$-annulus formula} is a formula $\phi$ of the form
$$
|p_0(x)| \Box_0 \epsilon_0 \wedge \bigwedge_{i=1}^L \epsilon_i \Box_i |p_i(x)|,
$$
where the $p_i \in K^\circ[x]$ are monic and irreducible, the
$\epsilon_i \in \sqrt{|K|\setminus \{0\}}$ and $\Box_i \in \{<,\le\}$.
Define $\overline\Box_i$ by $\{\Box_i, \overline\Box_i\}=\{<, \le\}$.
We require further that the \lq\lq holes'' $\{x\in \Kalg :
|p_i(x)| \overline\Box_i \epsilon_i\}$, $i=1,\dots,L$,
 all are contained in the
disc $\{x\in \Kalg:|p_0(x)|\Box_0 \epsilon_0\}$ and that the holes
corresponding to different indices $i$ are disjoint.
(Alternatively, we could require that $\epsilon_i \in
|K^\circ|\setminus \{0\}$ and allow the $p_i$ to be powers of
irreducible monic polynomials.)

(b) The corresponding {\em $K$-annulus} is
$$
\cU_\phi:= \{x \in \Kalg:\phi(x)\}
$$
(If $K_1 \supset \Kalg$ then $\phi$ also defines an  annulus in
$K_1$. We shall also refer to this
as $\cU_\phi$. No confusion will result.)

(c) a $K$-annulus formula $\phi$, and the corresponding
$K$-annulus $\cU_\phi$, is called {\em linear} if the $p_i$ are
all linear. (If $K=\Kalg$ then all $K$-annulus formulas are
linear.)

(d) a $K$-annulus fomula $\phi$, and the corresponding $K$-annulus
$\cU_\phi$, is called \emph{closed} (respectively \emph{open}) if
all the $\Box_i$ are $\le$ (respectively $<$).
\end{defn}

\begin{lem}
\label{annulus_properties}
(i) Let $p\in K[x]$ be irreducible and let $\Box\in\{<,\le\}$. Then for every $\delta\in\sqrt{\abs{K\setminus\{0\}}}$ there is an $\epsilon\in\sqrt{\abs{K\setminus\{0\}}}$ such that for every $x\in\Kalg$, $\abs{p(x)}\Box\varepsilon$ if, and only if, for some zero $\alpha$ of $p$, $\abs{x-\alpha}\Box\delta$.

(ii) A $K$-annulus is a finite union of isomophic (and linear) $\Kalg$-annuli.

(iii) Any two $K$-discs (i.e. $L=0$ in Definition~\ref{Rational
definitionA}) $\cU_1$ and $\cU_2$ are either disjoint or one is
contained in the other.

(iv) For any two $K$-annuli $\cU_1$ and $\cU_2$, if $\cU_1
\cap \cU_2 \not=\emptyset$
then $\cU_1 \cap \cU_2$ is a $K$-annulus.

(v) The complement of a $K$-annulus is a finite union of $K$-annuli.

(vi) Every set of the form
$$
\cU=\left\{x\in\Kalgo: |p_0(x)|\, \Box_0\, \epsilon_0 \wedge
\bigwedge^s_{i=1}\,\epsilon_i\,\Box_i|p_i(x)|\right\}
$$
with the $p_i$ irreducible is described by a $K$-annulus formula.
\end{lem}

\begin{proof}
Exercise
\end{proof}

\begin{defn}\label{def:holo}Let $\phi$ be a $K$-annulus formula as
in Definition \ref{Rational definitionA} (a). Define the {\em ring
of $K$-valued functions $\cO_K(\phi)$ on} $\cU_\phi$ by
$$
\cO_K(\phi):= S_{m+1,n}(\sigma,K)/ (p_0^{l_0}(x) - a_0z_0,
 p_1^{l_1}(x)z_1 -
a_1,\dots,p_L^{l_L}(x) z_L - a_L),
$$
where $a_i \in K^\circ$, $|a_i|= \epsilon_i^{l_i}$, $m+n=L+1$,
$\{z_0, \dots, z_L\}$  is the set $\{\xi_2, \dots, \xi_{m+1},
\rho_1, \dots, \rho_n\}$ and $x$ is $\xi_1$ and $z_i$ is a $\xi$ or
$\rho$ variable depending, respectively, on whether $\Box_i$ is
$\le$ or $<$. Observe that each $f\in \cO_K(\phi)$ defines a
function $\cU_\phi \to \Kalg$ via the analytic structure on $K$,
which by Theorem \ref{algebraic_extension} extends uniquely to
$\Kalg$.
\end{defn}

\begin{rem}
\label{cor:hookright} Let $\phi$ be a $K$-annulus formula and
suppose $\cU_\phi=\dot{\bigcup}\,\cU_i$, where the $\cU_i$ are the
(linear) $\Kalg$-annuli that comprise the  irreducible components
of $\cU_\phi$ over $\Kalg$ (cf.
Lemma~\ref{annulus_properties}~(ii)). Then one can prove
$$
\cO_K(\phi)\hookrightarrow \Kalg\otimes_K\cO_K(\phi)=\bigoplus_i
\cO_{\Kalg}(\cU_i).
$$
This result is not needed here.
\end{rem}

\begin{defn}\label{strongunit}
Let $f$ be a unit in $\cO_K(\phi)$. Suppose that there is some
$\ell\in\mathbb{N}$ and $c\in K$ such that $\abs{f^\ell(x)}=\abs c$
for all $x\in\cU_\phi$. Suppose also that there exists a polynomial
$P(\xi)\in\Kt[\xi]$ such that
$P\left(\left(\frac{1}{c}f^\ell(x)\right)^\sim\right)=0$ for all
$x\in\cU_\phi$, where $^\sim:\Ko\to \Kt$ is the natural projection.
Then call $f$ a \emph{strong unit}. Call $f$ a \emph{very strong}
unit if moreover $\abs{f(x)}=1$ and $\left(f(x)\right)^\sim=1$ in
$\cO_K(\phi)$.
\end{defn}

\begin{lem}[Normalization]
\label{lem:normalization}
(i) Let $\phi$ be a closed $K$-annulus formula. Then there is an inclusion
$$
S_{1,0}(\sigma,K) \hookrightarrow \cO_K(\phi),
$$
which is a finite ring extension.

(ii) Let $\phi$ be an open $K$-annulus formula. Then there is an
inclusion
$$
S_{0,1}(\sigma,K) \hookrightarrow \cO_K(\phi),
$$
which is a finite ring extension.
\end{lem}

\begin{proof}
Apply a suitable Weierstrass automorphism, as in the classical case.
\end{proof}

The following two corollaries are proved exactly as in the
classical case (cf. \cite{bgr}).

\begin{cor}
\label{cor:theNullstellensatz}
Let $\phi$ be a $K$-annulus
formula that is either closed or open. Then

(i) the Nullstellensatz holds for $\cO_K(\phi)$; i.e., the maximal ideals of $\cO_K(\phi)$ are $K$-algebraic.

(ii) $\cO_K(\phi)$ is an integral domain.
\end{cor}

\begin{cor}
\label{cor:Nulle}
(i) If $\phi$ is a $K$-annulus formula that is either closed or open and
$\psi$ is any $K$-annulus formula with $\cU_\phi \subseteq \cU_\psi$ then
$$
\cO_K(\psi)\hookrightarrow \cO_K(\phi) = \cO_K(\cU_\phi).
$$

(ii) If $\phi$ is a $K$-annulus formula that is closed or open, the ring $\cO_K(\phi)$
depends only on $\cU_\phi$ and is independent of the presentation
$\phi$.
\end{cor}

Let $U$ be a $K$-annulus. A $K$-annulus formula $\phi$
$$
|p_0(x)| \Box_0 \epsilon_0 \wedge \bigwedge_{i=1}^L \epsilon_i \Box_i |p_i(x)|
$$
is called a {\em good description} of $U$ if, and only if,
$U=\cU_\phi$ and each $p_i$ is of minimal degree. This condition
implies that if $\deg q<\deg p_i$ then $q$ has no zero in the hole
defined by $p_i$; i.e., in the disc defined by the formula
$\abs{p_i(x)}\overline\Box_i\epsilon_i$. If $\phi$ is a good
description of $\cU_\phi$, then we say that $\phi$ is a {\em good
$K$-annulus formula}. Observe that each $K$-annulus has a good
description. Moreover, by Corollary~\ref{cor:Nulle}~(ii), if
$\phi$ is a closed or open $K$-annulus formula, then replacing
$\phi$ by a good description does not change the ring of analytic
functions.

\par
The main result of this section is the following.
\begin{thm}
\label{thm:torational} Let $\phi$ be a $K$-annulus formula and let
$f\in\cO_K(\phi)$. Then there are: finitely many $K$-annulus
formulas $\phi_i$, each either closed or open, such that $\phi$ is
equivalent to the disjunction of the $\phi_i$, rational functions
$R_i\in\cO_K(\phi_i)$ and very strong units $E_i\in\cO_K(\phi_i)$
such that for each $i$,$f|_{\cU_{\phi_i}}=R_iE_i$.
\end{thm}

\begin{proof}
By Lemma~\ref{annulus-decomposition}, the theorem follows from
Propositions~\ref{thin-factorization}
and~\ref{Laurent-factorization}, below.
\end{proof}

The decomposition in Theorem~\ref{thm:torational} will be given in
terms of two types of annuli, thin annuli and Laurent annuli,
defined as follows.

\begin{defn}
\label{thin-Laurent}
(i) A linear $K$-annulus is called {\em thin} if, and only if, it is of the form
$$
\{x\in\Kalg:\abs{x-a_0}\le\epsilon\ \mathrm{and\ for}\ i=1,\dots,n,\ \abs{x-a_i}\ge\epsilon\}
$$
for some $\epsilon\in\sqrt{\abs{K\setminus\{0\}}}$,
$\epsilon\le1$, and $a_i\in\Ko$. Recall, by
Lemma~\ref{annulus_properties}, that a $K$-annulus $\cU$ is a
union of isomorphic (and linear) $\Kalg$-annuli $\cU_i$. In
general, the $K$-annulus $\cU$ is called {\em thin} if, and only
if, each $\cU_i$ is a thin $\Kalg$-annulus.

(ii) A $K$-annulus $\cU$ of the form
$$
\{x\in\Kalg:\epsilon_1<\abs{p(x)}<\epsilon_0\},
$$
where $p\in K[x]$ is irreducible is called a {\em Laurent annulus}.
\end{defn}

Note that any open $K$-annulus with one hole is Laurent.

\begin{lem}
\label{annulus-decomposition} Every $K$-annulus  formula is
equivalent to a finite disjunction of $K$-annulus formulas that are
either thin or Laurent.
\end{lem}

\begin{proof}
Define the {\em complexity} of a $K$-annulus formula
$$
\abs{p_0(x)}\Box_0\epsilon_0\wedge\bigwedge_{i=1}^L\epsilon_i\Box_i\abs{p_i(x)}
$$
to be $\sum\limits_{i=1}^L\deg p_i$. The proof is by induction on complexity. The base case is easy.

By removing thin annuli, we may assume that the remaining set is an
open $K$-annulus. If there is only one hole, the annulus is Laurent,
and we are done. Assume that there are at least two holes. Two
closed holes are said to {\em abut} when their radii are equal to
the distance between the centers. After removing the largest
possible Laurent annuli surrounding each hole, we may assume that
all holes are closed, and at least two abut. Now, removing a  thin
annulus lowers the complexity.
\end{proof}

\begin{lem}
\label{sum-max}
Consider the following $K$-annulus formula, $\phi$:
$$
\abs{p_0(x)}\le\epsilon_0\wedge\bigwedge_{i=1}^L\epsilon_i\le\abs{p_i(x)}.
$$
Suppose that $\phi$ is a good, thin $K$-annulus formula. Let
$\nu_{ij}\in\mathbb{N}$, $i=1,2$, $j=1,\dots,L$. Suppose $f_i\in
K[x]$ satisfy $\deg f_i<\deg p_j$ for all $j$ such that $\nu_{ij}>0$
and suppose that $\nu_{1j}\not=\nu_{2j}$ for some $j$. Then
$$
\norm{\frac{f_1}{\prod_j p_j^{\nu_{1j}}}+\frac{f_2}{\prod_j p_j^{\nu_{2j}}}}=
\max\left\{\norm{\frac{f_1}{\prod_j p_j^{\nu_{1j}}}},\norm{\frac{f_2}{\prod_j p_j^{\nu_{2j}}}}\right\}.
$$
\end{lem}

\begin{proof}
This reduces easily to the linear case, which is treated in the
proof of \cite{fp}, Theorem 2.2.6.
\end{proof}

\begin{lem}
\label{sup-valuation} When $\phi$ is a  thin $K$-annulus formula,
$\norm\cdot$ is a valuation on $\cO_K(\phi)$.
\end{lem}

\begin{proof}
Let the notation be as in the statement of Lemma~\ref{sum-max}.
By that lemma and the definition of $\cO_K(\phi)$, this reduces to showing that
$$
\norm{\frac{f_1}{\prod p_j^{\nu_{1j}}}\cdot\frac{f_2}{\prod
p_j^{\nu_{2j}}}}= \norm{\frac{f_1}{\prod
p_j^{\nu_{1j}}}}\cdot\norm{\frac{f_2}{\prod p_j^{\nu_{2j}}}},
$$
which is immediate.
\end{proof}

\begin{prop}
\label{thin-factorization} Let $\phi$ be a thin $K$-annulus formula.
Then for each $f\in\cO_K(\phi)$ there is a rational function
$R\in\cO_K(\phi)$ and a very strong unit $E\in\cO_K(\phi)$ such that
$f=R\cdot E$.
\end{prop}

\begin{proof}
By the Nullstellensatz, there is a monic polynomial $f_0\in\Ko[x]$,
with zeros only in $\cU_\phi$, and an $f'\in\cO_K(\phi)$ such that
$f=f_0\cdot f'$, and $f'$ is a unit of $\cO_K(\phi)$. Thus we may
assume that $f$ is a unit. We may also assume that $\norm f=1$. By
Lemma~\ref{sup-valuation}, this implies that $\norm g=1$, where
$g\in\cO_K(\phi)$ satisfies $gf=1$. Thus, by Lemma~\ref{sum-max} and
the definition of $\cO_K(\phi)$, there is a rational function $\hat
f\in\cO_K(\phi)$ such that $\norm{f-\hat f}<1$ and $\hat f$ is a
unit of $\cO_K(\phi)$. Since $\hat f$ is a unit, we may write
$f=\hat f\cdot E$. We have $\norm{E-1}<1$, so $E$ is a very strong
unit.
\end{proof}

\begin{prop}
\label{Laurent-factorization} Let $\phi$ be a Laurent $K$-annulus
formula, and let $f\in\cO_K(\phi)$. There are finitely many
$K$-annulus formulas $\phi_i$, each either thin or Laurent, such
that $\cU_\phi=\cup_i\cU_{\phi_i}$, and for each $i$, there are
rational functions $R_i\in\cO_K(\phi_i)$ and very strong units
$E_i\in\cO_K(\phi_i)$ such that $f|_{U_{\phi_i}}=R_i\cdot E_i$.
\end{prop}

\begin{proof}
Write
$$
f=\sum_{i\in\mathbb{Z}}a_i(x)p^i,
$$
where $p$ is the polynomial that occurs in $\phi$ and the $a_i$ are
polynomials of degree less than the degree of $p$. By
Lemma~\ref{strong_noetherian}, there are only finitely many $a_i$
that can be dominant (in the sense of the proof of
Proposition~\ref{kernel}) on any sub-annulus of $\cU_\phi$. There is
a partition  of $\cU_\phi$ into a finite collection of thin
sub-annuli and Laurent sub-annuli such that each Laurent sub-annulus
is either of lower complexity or on the Laurent sub-annulus, each of
the finitely many dominant $a_i$ is a strong unit (only the fact
that it is a unit is used). The thin sub-annuli are handled by
Proposition~\ref{thin-factorization}. The Laurent sub-annuli of
lower complexity are treated by induction and the remaining Laurent
sub-annuli are treated as in \cite{lp}, Theorem~3.3.
\end{proof}

\begin{lem}
\label{D-function}
Let
$$
R(x)=x^{n_0}\,\prod^s_{i=1}\, p_i(x)^{n_i} \in K(x),
$$
where the $p_i\in \Ko[x]$ are monic, irreducible and mutually
prime and the $n_i\in \ZZ$. Let $\epsilon \in
\sqrt{|K|\setminus\{0\}}$, let $\Box \in \{<,\le\}$, and let
$$
\cU:=\{x\in\Kalgo : |R(x)|\, \Box \, \epsilon\}.
$$
There are finitely many $K$-annuli $\cU_i$, $i=1,\dots,L$, such that
$\cU=\bigcup^L_{i=1}\cU_i$ and for each $i$,
$R(x)|_{\cU_i} \in \cO_K(\cU_i)$.
\end{lem}

\begin{proof}
Induction on $s$. For each $i$, let $\alpha_i \in \Kalgo$ be a zero of
$p_i$, and let $a_i:=|\alpha_i|\in \sqrt{|K|}$. Since $p_i$
is irreducible, $a_i$ is independent of which zero of $p_i$ is chosen. Let
$$
d_{ij} := \min \{|\alpha-\beta|: p_i(\alpha)=0=p_j(\beta)\}\in \sqrt{|K|}.
$$
Hence $d_{ij}$ is the smallest distance between a zero of $p_i$ and a
zero of $p_j$. We consider several cases.

\item[(Case 1)] There are $i,j$ with $a_i<a_j$ . Choose $\gamma \in
\sqrt{|K|}$ with $a_i<\gamma <a_j$. Let $\cU_1 :=
\cU\cap\{x:|x|\le\gamma\}$ and $\cU_2 := \cU\cap\{x:|x|\ge\gamma\}$.
Then on $\cU_1$, $p_j$ is a strong unit and on $\cU_2$, we have
that $p_i=E\cdot x^{deg(p_i)}$, $E$ a strong unit.
\item[(Case 2)] $a_i=\gamma$ for all $i$ and $d_{12}< d_{13}$, say.
Choose $\delta\in \sqrt{|K|}$ so that $d_{12}<\delta< d_{13}$ and by
Lemma~\ref{annulus_properties}~(i), choose $\gamma'\in \sqrt{|K|}$ so that
$$
|p_1(x)|\le \gamma' \longleftrightarrow \bigvee_{\scriptstyle{\alpha
\text{ such that} \atop \scriptstyle p_1(\alpha)=0}} \, |x-\alpha|\le
\delta.
$$
On $\cU_3:= \cU \cap \{x:|p_1(x)|\le \gamma'\}$, $p_3$ is a
strong unit. On $\cU_4:= \cU \cap \{x:\gamma'\le |p_1(x)|\}$, we
have that $p_2(x)^{n_1}=E\cdot p_1(x)^{n_2}$, $E$ a strong unit,
for suitable $n_1, n_2 \in \mathbb{N}$.

\item[(Case 3)] $a_i=\gamma$ for all $i$ and $d_{ij}=\delta$ for all $i\ne j$.
Then $\delta \le \gamma$. Choose $\gamma_i$ such that
$$
\{x:|p_i(x)|<\gamma_i\} = \bigcup_{\scriptstyle{\alpha\text{ such that }} \atop
\scriptstyle p_i(\alpha)=0} \, \{x:|x-\alpha|< \delta\}
$$
On $\cU_{5i}:=\{x:|p_i(x)|< \gamma_i\}$, each $p_j$ with $j \ne i$ is
a strong unit. On $\cU_{6i}:=\{x:|p_i(x)|= \gamma_i\}$, each $p_i$ is of constant
size $\gamma_i$. On $\cU_7:=\{x:|p_1(x)|> \gamma_1\}$  there are $d_i \in \QQ_+$
such that $|p_i(x)|=|p_1(x)|^{d_i}$. (For each zero $\alpha_i$
of $p_i$ there is a zero $\alpha_1$ of $p_1$ so that for all $x \in
\cU_7$ we have $|x-\alpha_i|=|x-\alpha_1|$.)
\end{proof}

Using Theorem~\ref{thm:torational}, Propositions
\ref{thin-factorization} and \ref{Laurent-factorization},
Lemma~\ref{D-function}, and induction on the complexity of terms, we
obtain the following.

\begin{thm}
\label{lem:toholom}
 Let $\tau(x)$ be an $\Lan$-term of the valued
field sort (cf.~section \ref{sec:language}). Then there is a finite
set $S\subset \Kalgo$ and a finite collection of disjoint $K$-annuli
$\cU_i$, each open or closed, and for each $i$, an $F_i\in
\cO_K(\cU_i)$ such that $\Kalgo = \bigcup \, \cU_i$ and
$$ \tau\bigm|_{\cU_i\setminus S}= F_i\bigm|_{\cU_i\setminus S}.$$
Moreover, we can ensure that
$$
F_i=R_iE_i
$$
where $R_i$ is a rational function over $K$ and $E_i\in
\cO_K(\cU_i)$ is a very strong unit.
\end{thm}

\begin{cor}
\label{cor:quantifier} Let $\psi(x)$ be a quantifier free $\LL_{\rm
Fields}\cup S_{m,n}(A)$-formula, where $\LL_{\rm Fields}$ is the
language of fields $(+,-,^{-1},\cdot,0,1)$. Then there is a
polynomial $F$ and (closed or open) $K$-annulus formulas $\phi_i$
such that for every field $K_1 \supset K$ with $A$-analytic
structure
$$K^\circ_1 \models F(x)=0 \vee \left[\psi(x) \longleftrightarrow
\bigvee_i \, \phi_i(x)\right]
$$
\end{cor}

\begin{cor}\label{cor:strictrational}
Let $K$ be a valued field whose value group has a minimal positive
element. Suppose that $K$ has a strictly convergent $A$-analytic
structure $\sigma$ as in Definition \ref{analytic_structure}.
Suppose that a uniformizer of $\Ko$ is given by $\sigma_0(v)$ for
some $v$ in $A$. Let $K_{\rm unram}$ be the maximal unramified
(algebraic) field extension of $K$.

Let $\tau(x)$ be an $\Lans$-term of the valued field sort
(cf.~section \ref{sec:language}) and let $n>0$ be an integer. Then
there is a finite set $S\subset K_{\rm unram}^\circ$ and a finite
collection of disjoint closed $K$-annuli $\cU_i$, and for each $i$,
an $F_i\in \cO_K(\cU_i)$, such that $K_{\rm unram}^\circ = \bigcup
\, \cU_i'$, where $\cU_i'=K_{\rm unram}\cap \cU_i$ and
$$ \tau\bigm|_{\cU'_i\setminus S}= F_i\bigm|_{\cU'_i\setminus S}.$$
Moreover, we can ensure that
$$
F_i=R_iE_i
$$
where $R_i$ is a rational function over $K$ and $E_i\in
\cO_K(\cU_i)$ is a very strong unit which furthermore satisfies
$$
E_i\equiv 1 \bmod \sigma_0(v)^n.
$$\end{cor}
\begin{proof}
Since we work in $K_{\rm unram}$, we can replace strict inequalities
by weak inequalities using the uniformizer $\sigma_0(v)$.
\end{proof}
\begin{rem}\label{rem:rat:strict}
Note that the decomposition of $K_{\rm unram}$ and the
representation of $\tau$ given Corollary \ref{cor:strictrational}
actually holds in all unramified field extensions of $K$ with
$A$-analytic structure extending $\sigma$.
\end{rem}

\section{$A$-analytic languages and quantifier elimination}\label{sec:language}

In this section we recall the notion of languages of Denef-Pas and
we introduce the notion of $A$-analytic languages, suitable for
talking about valued fields with $A$-analytic structure. Further, we
specify the theories that we will consider, and we establish the
corresponding quantifier elimination results in equicharacteristic
zero and in mixed characteristic with bounded ramification.

\par
For $K$ a valued field, $I$ an ideal of $K^\circ$, write ${\rm
res}_I:K^\circ\to K^\circ/I$ for the natural projection. An
\emph{angular component modulo $I$} is a map $\ac_I:K\to K^\circ/I$
such that the restriction to $K^\times$ is a multiplicative
homomorphism to $(K^\circ/I)^\times$, the restriction to
$(K^\circ)^\times$ coincides with the restriction to
$(K^\circ)^\times$ of ${\rm res}_I$, and such that $\ac_I(0)=0$.

\par
Fix a sequence of positive numbers $(n_p)_{p}$, indexed by the prime
numbers and write $\NN_0:=\{x\in\ZZ :  x>0\}$. We consider
structures
$$(K,\{K^\circ/I_m\}_{m\in \NN_0},\ord(K^\times)),$$
where $K$ is a Henselian valued field of characteristic zero with
valuation ring $\Ko$, additively written valuation\footnote{Our use
of the symbol $\ord$ with argument $x$ in an $\LPas$-formula is in
fact an abbreviation for a function with domain the ${\rm Val}$-sort
which extends the valuation (the reader may choose the value of
$0$), conjoined with the condition $x\not=0$.}
$\ord:K^\times\to\ord(K^\times)$, angular component maps $\ac_{m}$
modulo $I_m$, a constant $t_K\in K$, and ideals $I_m$ of $K^\circ$
for $m\in \NN_0$, satisfying the following properties
\begin{itemize}
 \item[(I)] $I_1$ is the maximal ideal of $\Ko$,
 \item[(II)] either $I_m=I_1$ for all $m\in\NN_0$ and $t_K=1$, or, $\ord(K^\times)$ has a
 minimal positive element, $I_m=I_1^m$ for all $m\in\NN_0$, and
 $t_K$ is either $1$ or an element of $\Ko$ with minimal positive valuation
 such that $\ac_m(t_K)=1$ for all $m>0$,
 \item[(III)]
 if the residue field $\Kt$ of $K$ has characteristic
$p>0$, then,  $t_K\not=1$ (hence, $I_2\not=I_1$), and the
ramification is bounded by $\ord(p)\leq n_p$.
\end{itemize}

\par
Let $\cK((n_p)_p)$ be the class of these structures. We call the
sorts $\Val$ for the valued field sort, $\Res_m$ for the $m$-th
residue ring $\Ko/I_m$ for $m\in \NN_0$, more generally $\Res$ for
the disjoint union of the $\Res_m$, and $\Ord$ for the value group
sort.

\par
For $\Val$ we use the language
$\LL_{\Val}=(+,-,^{-1},\cdot,0,1,t_0)$ of fields with an extra
constant symbol $t_0$, interpreted in $K$ as $t_K$. Let $\LL_{{\rm
Ord},0}=(+,-,\leq,0)$ be the language of ordered groups.

\par
For $\Res$ we define the language $\LL_{{\rm Res},0}$ as the
language having the ring language and a constant symbol $t_m$ for
each sort $\Res_m$ and natural projection maps $\pi_{mn}:\Res_m \to
\Res_n$ giving commutative diagrams with the maps ${\rm res}_m:={\rm
res}_{I_m}$ and ${\rm res}_n$ for $m\geq n$. If $I_2=I_1$, the $t_m$
are interpreted as $1$. If $I_2\not=I_1$, $t_m$ denotes the image
under ${\rm res}_m$ of an element $x$ with $\ac_m(x)=1$ and
$\ord(x)$ the minimal positive element.

\par
Fix expansions $\LL_{\rm Ord}$ of $\LL_{{\rm Ord},0}$ and $\LL_{{\rm
Res}}$ of $\LL_{{\rm Res},0}$.
 To this data we associate the \emph{language
$\LPas=\LPas(\LL_{{\rm Res}},\LL_{\rm Ord})$ of Denef-Pas} defined
as
 $$(\LL_{\Val},\LL_{{\rm Res}},\LL_{\rm Ord},
 \{\ac_m\}_{m\in \NN_0},\ord).$$

\par
Fix a language $\LPas$ of Denef-Pas, an $\LL_{\Ord}$-theory
$\mathbf{T}_{\Ord}$, and an $\LL_{\Res}$-theory $\mathbf{T}_{\Res}$.
 For such data, let $\TPas=\TPas(\LPas,\mathbf{T}_{\Res},\mathbf{T}_{\Ord},(n_p)_{p})$
be the $\LPas$-theory of all structures
$$(K,\{K^\circ/I_m\}_{m\in \NN_0},\ord(K^\times))$$
in $\cK((n_p)_{p})$ which are $\LPas$-structures such that
$\{K^\circ/I_m\}_{m\in \NN_0}$ is a model of
$\mathbf{T}_{{\Res}}$, and $\ord(K^\times)$ is a model of
$\mathbf{T}_{\rm Ord}$.

\par
The sorts $\Res_m$ for $m\in \NN_0$ and $\Ord$ are called
\emph{auxiliary sorts}, and $\Val$ is the main sort.

\par
Now we come to the notion of $A$-analytic languages. Fix a
Noetherian ring $A$ that is complete and separated for the
$I$-adic topology for a fixed ideal $I$ of $A$. Define the
\emph{separated $A$-analytic language}
$$\Lan:=\LPas\cup_{m,n\geq 0} S_{m,n}(A) $$
and the \emph{strictly convergent $A$-analytic language}
$$\Lans:=\LPas\cup_{m\geq 0}
T_{m}(A),
$$
where $S_{m,n}(A)$ and $T_m(A)$ are as in Definition
\ref{strictly_convergent}.

\par
For $\TPas$ as before, we define the $\Lan$-theory
$$\Tan:=\TPas\cup \mathrm{(IV)}_{S}$$
and the $\cL_{T(A)}$-theory
$$\Tans:=\TPas\cup \mathrm{(IV)}_{T},$$
where
\begin{itemize}
 \item[(IV)$_{S}$]
  the $\Val$-sort is equipped with a separated analytic
 $A$-structure and each symbol $f$ of $S_{m,n}(A)$ is considered as
 a function $\Val^{m+n}\to \Val$ by extending $f$ by zero
 outside its domain $\{(x,y)\in \Val^{m+n} :  \ord(x_i)\geq 0,\ \ord(y_j)>0\}$,
 \item[(IV)$_{T}$]
  the $\Val$-sort is equipped with a strictly convergent analytic
 $A$-structure and each symbol $f$ of $T_{m}(A)$ is considered as
 a function $\Val^{m}\to \Val$ by extending $f$ by zero
 outside its domain $\{x\in \Val^m :  \ord(x_i)\geq 0\}$.
 Moreover, the value group has a minimal
 positive element and this is the order of some constant of
 $T_0(A)$.
\end{itemize}

Observe that there exist $\Lan$-terms, resp.~$\Lans$-terms, which
yield all kinds of restricted division as considered in \cite{mod}
and \cite{unifrid}.

\par
Write ${\rm res}_m$ for the natural projection from the valuation
ring to $\Res_m$, for each $ m>0$. Note that each map ${\rm res}_m$
is definable without $\Val$-quantifiers, since it sends $x\in \Val$
to $\ac_m(x)$ when $\ord(x)=0$, and to $\ac_m(1+x)-1$ when
$\ord(x)>0$.

\par
For each language $\cL$, write $\cL^*$ for the expansion
\begin{equation}\label{Lstar}
\cL^*:=\cL\cup_{m>0,\ e\geq 0}\{(\cdot,\cdot,\cdot)_e^{1/m},\
h_{m,e}\},
 \end{equation}
 with $m,e$ integers.
 Then, each model of $\TPas$ extends uniquely to an
$\LPas^*$-structure, axiomatized as follows:
\begin{itemize}
  \item[(V)]
   $(\cdot,\cdot,\cdot)_e^{1/m}$
 is the function
  $\Val\times \Res_{2e+1}\times \Ord \to \Val$
  sending $(x,\xi,z)$ to the (unique) $m$-th root $y$ of
  $x$ with $\ac_{e +1}(y)=\pi_{2e +1,e+1}(\xi)$ and $\ord(y)=z$, whenever
  $\xi^m=\ac_{2e +1}(x)$,
  $m\not=0$ in $\Res_{e +1}$,
  and $mz=\ord(x)$, and to $0$
  otherwise;
  \item[(VI)]
   $h_{m,e}$ is the function
  $\Val^{m+1}\times \Res_{2e +1} \to \Val$
sending $(a_0,\ldots,a_{m},\xi)$ to the unique $y$ satisfying
$\ord(y)=0$, $\ac_{e +1}(y)=\pi_{2e +1,e +1}(\xi)$, and
$\sum_{i=0}^{m} a_{i} y^i=0$, whenever $\xi$ is a unit,
$\ord(a_i)\geq 0$, $f(\xi)=0$, and
 $$\pi_{2e +1,e +1}(f'(\xi))\not=0,$$
with $f(\eta)=\sum_{i=0}^m{\rm res}_{2e+1}(a_{i}) \eta^i$ and $f'$
its derivative, and to $0$ otherwise.
\end{itemize}

\par
Sometimes we will use the property, for $\ell$ either zero or a
prime number,
\begin{itemize}
 \item[(VII)$_\ell$]
the residue field has characteristic $\ell$.
\end{itemize}

\par
The following result extends quantifier elimination results of van
den Dries \cite{vdd1} and Pas in \cite{pas} and \cite{pas2}.
 Theorem \ref{thm:qe:an} for the theory $\TPas$ can be compared with
results obtained by Kuhlmann \cite{Kuhl}. In \cite{Kuhl}, the
language for the auxiliary sorts is less explicit than in this
paper.

\begin{thm}[Quantifier elimination]\label{thm:qe:an}
Let $\mathcal{T}$ be one of the theories $\TPas$, $\Tan$, or
$\Tans$, and let $\cL_{\cT}$ be respectively $\LPas$, $\Lan$, or
$\Lans$. Then $\mathcal{T}$ admits elimination of quantifiers of the
valued field sort. Moreover, every $\cL_{\cT}$-formula $\phi (x,
\xi, \alpha)$, with $x$ variables of the valued field-sort, $\xi$
variables of the residue rings sorts and $\alpha$ variables of the
value group sort, is $\mathcal{T}$-equivalent to a finite
disjunction of formulas of the form
\begin{equation}\label{form}
\psi (\ac_\ell f_1 (x), \dots, \ac_\ell f_k (x), \xi) \wedge
\theta (\ord\, f_1 (x), \dots, \ord\, f_k (x), \alpha),
\end{equation}
with $\psi$ an $\LL_{\Res}$-formula, $\theta$ an $\LL_{\rm
Ord}$-formula, and $f_1 \dots, f_k$ $\cL_{\cT}$-terms.
\end{thm}
\begin{proof}[Proof of Theorem \ref{thm:qe:an}]
If one knows the quantifier elimination statement, the statement
about the form of the formulas follows easily, cf.~\cite{vdd1}.

\par
The quantifier elimination statement for $\TPas$ is proved together
with the cell decomposition Theorem \ref{thm:cd} for $\TPas$.
\par
The statement for the analytic theories follows in a nowadays
standard way from the result for $\TPas$ and the Weierstrass
division as developed in Section~\ref{part:an}, cf.~the proof of
Theorem 3.9 of \cite{vdd1} in the strictly convergent case, and
Theorem 4.2 of \cite{unifrid} in the separated case.
\end{proof}

\section{Definable assignments}\label{sec:assignm}

We elaborate on the terminology of \cite{DL} and \cite{CLoes} on
definable assignments and definable subassignments.  By some
authors, definable assignments are just called ``formulas'', or
``definable sets'', and definable morphisms are often called
``definable functions''.

\par
Let $\cT$ be a multisorted theory formulated in some language $\cL$,
where some of the sorts are auxiliary sorts and the other sorts are
main sorts. Let $\Mod(\cT)$ be the category whose objects are models
of $\cT$ and whose morphisms are elementary embeddings. By a
\emph{$\cT$-assignment} $X$ we mean a $\cT$-equivalence class of
$\cL$-formulas $\varphi$, where we say that $\varphi$ and $\varphi'$
are $\cT$-equivalent if they have the same set of $\cM$-rational
points for each $\cM$ in $\Mod (\cT)$. Knowing a $\cT$-assignment
$X$ is equivalent to knowing the functor from $\Mod(\cT)$ to the
category of sets which sends $\cM\in\Mod(\cT)$ to the set
$\varphi(\cM)$ for any $\cL$-formula $\varphi$ in $X$, and we will
identify $\cT$-assignments with these functors.

\par
The usual set theoretic operations can be applied to
$\cT$-assignments, for example, for two $\cT$-assignments $X,Y$,
$X\subset Y$ has the natural meaning and if $X\subset Y$ call $X$ a
$\cT$-subassignment of $Y$. Similarly, for $X,Y\subset Z$
$\cT$-assignments, $X\cup Y$, $X\cap Y$, and $X\setminus Y$ have the
obvious meaning. Cartesian products have the obvious meaning and
notation. We refer to \cite{CLoes} and \cite{DL} for more details on
the general theory of assignments and $\cT$-assignments.

\par
For $\cT$-assignments $X,Y$, a collection of functions
$f_{\cM}:X(\cM)\to Y(\cM)$ for each $\cM\in\Mod(\cT)$ is called a
\emph{$\cT$-morphism} from $X$ to $Y$ if the functor sending
$\cM\in\Mod(\cT)$ to the graph of $f_{\cM}$ is a $\cT$-assignment. A
$\cT$-morphism $f:X\to Y$ such that $f_{\cM}$ is a bijection for
each $\cM\in\Mod(\cT)$ is called a $\cT$-isomorphism.

\begin{defn}\label{def:parameterization}
By a $\cT$-\emph{parameterization} of a $\cT$-assignment $X$, we
mean a $\cT$-isomorphism $f:X\to Y\subset X\times R$ with $R$ a
Cartesian product of auxiliary sorts, such that $\pi\circ f:X\to
X$ is the identity on $X$, with $\pi$ the projection.
\end{defn}

\begin{ex}
If $\cT$ is one of the theories $\TPas$, $\Tan$, or $\Tans$, write
$\Val^{\ell_1}\times \Res_n^{\ell_2} \times \Ord^{\ell_3}$, for the
$\cT$-assignment which sends a model
$$(K,\{K^\circ/I_m\}_{m>0},\ord(K^\times))$$
  to
$$K^{\ell_1}\times
(K^\circ/I_n)^{\ell_2}\times\ord(K^\times)^{\ell_3},$$
 for any
$n>0,\ell_i\geq 0$. We recall that, for such $\cT$, the sorts
$\Res_m$ and $\Ord$ are called auxiliary sorts. For such $\cT$, the
map $\Val\to\Val\times\Res_1:x\mapsto (x,\ac_1(x))$ is an example of
a $\cT$-parameterization.
\end{ex}

\section{Cell decomposition}\label{sec5}

In this section we state and prove an analytic cell decomposition
theorem for $\cT$-assignments with $\cT$ one of the theories
$\Tan$ or $\Tans$, see Theorem \ref{thm:cd} below.
 Theorem \ref{thm:cd} generalizes cell decompositions of \cite{pas},
\cite{pas2}, \cite{Denef2}, \cite{c}, and \cite{CLoes} and provides
what is needed for the applications to analytic integrals in the
next section.
 Also for the theory $\TPas$ we obtain a cell decomposition,
cf.~Theorem \ref{thm:cd}, which generalizes and refines the cell
decompositions of Pas \cite{pas}, \cite{pas2} in several ways: in
equicharacteristic zero, also angular components of higher order
(i.e., modulo powers of the maximal ideal) are allowed; we can take
the centers of the cells to be $\LPas^*$-terms; we can partition any
definable set into cells adapted to any given definable
function\footnote{The cells used by Pas are neither suitable for the
partition of definable sets, nor for the preparation of definable
functions. This problem has been addressed in \cite{CLoes}.}; in
mixed characteristic, we allow for any value group with a least
positive element\footnote{In mixed characteristic, Pas \cite{pas2}
allows for the integers as value group only.}. In \cite{CLoes}, a
notion of cells is introduced that is more general than the one in
\cite{pas}; we base the definition of cells below on this notion of
\cite{CLoes}. Summarizing, our cell decomposition holds for the
theories introduced in section \ref{sec:language}, which includes
both the analytic and the algebraic cases, in equicharacteristic
zero with quite general angular components, and in mixed
characteristic as long as the degree of ramification is bounded.
Crucial in the proof of analytic cell decomposition are Propositions
\ref{thin-factorization} and \ref{Laurent-factorization}, Theorem
\ref{lem:toholom}, and Corollary \ref{cor:strictrational}. The proof
of the analytic cell decomposition seems to require all the work of
the previous sections. As a second main result of this section, we
prove the fundamental result that $\cL$-definable functions are,
after parameterization using auxiliary sorts, given by
$\cL^*$-terms, where $\cL$ is either $\LPas$, $\Lan$, or $\Lans$,
see Theorem \ref{normal:an}; this result is motivated by notes of
van den Dries \cite{vdd2}.

\par
Fix $\cT$ to be $\TPas$, $\Tan$, or $\Tans$, and let $\cL_{\cT}$ be
respectively $\LPas$, $\Lan$, or $\Lans$. We come to the somewhat
elaborate definition of $\cT$-cells. (More generally, for such
$\cT$, the following definition makes sense as well for the
$\cL_\cT\cup K$-theory $\cT(K)$, where $K$ is the valued field of a
model of $\cT$, and $\cT(K)$ the theory $\cT$ together with the
diagram of $K$.)
\begin{defn}\label{defcells}
Let $C$ be a $\cT$-assignment, $k>0$ an integer, and $\alpha : C
\rightarrow \Ord$, $\xi : C \rightarrow \Res_k$, and $c : C
\rightarrow \Val$ $\cT$-morphisms, such that $\xi$ is always a
multiplicative unit in the $\Res_k$-sort. The $1$-$\cT$-cell
$Z_{C, \alpha, \xi, c}$ with base $C$, order $\alpha$, angular
component $\xi$, and center $c$ is the $\cT$-subassignment of
$C\times\Val$ defined by
   $$
y \in C\ \wedge\ \ord\, (z - c (y)) = \alpha (y)\ \wedge\ \ac_k (z
- c(y)) = \xi (y),
   $$
where $y\in C$ and $z\in \Val$. Similarly, if $c$ is a
$\cT$-morphism $c : C \rightarrow \Val$, we define the
$0$-$\cT$-cell $Z_{C, c}\subset C\times \Val$ with base $C$ and
center $c$ as the $\cT$-subassignment of $C\times\Val$ defined by
   $$
y \in C \wedge z = c(y).
   $$
More generally, $Z\subset S\times\Val$ with $Z$ and $S$
$\cT$-assignments will be called a $\cT$-1-cell, resp.~a
$\cT$-0-cell, if there exists a $\cT$-parameterization
$$
\lambda : Z \rightarrow Z_C \subset S\times R\times\Val,
$$
for some Cartesian product $R$ of auxiliary sorts and some
$\cT$-$1$-cell $Z_C=Z_{C, \alpha, \xi, c}$, resp.~$\cT$-$0$-cell
$Z_C=Z_{C, c}$.

We shall call the data $(\lambda, Z_{C, \alpha, \xi, c})$, resp.
$(\lambda, Z_{C, c})$, sometimes written for short $(\lambda, Z_C)$,
a \emph{$\cT$-presentation} of the $\cT$-cell $Z$.
\end{defn}

\begin{defn}
A $\cT$-morphism $f:Z\subset S\times\Val\to R$ with $Z$ a
$\cT$-cell, $S$ a $\cT$-assignment, and $R$ a Cartesian product of
auxiliary sorts, is called \emph{$\cT$-prepared} if there exist a
$\cT$-presentation $\lambda:Z\mapsto Z_C$ of $Z$ onto a cell $Z_C$
with base $C$ and a $\cT$-morphism $g:C\to R$ such that $f= g \circ
\pi \circ \lambda$, with $\pi:Z_C\to C$ the projection.
\end{defn}

\begin{ex}
Let $Z$ be the $\cT$-cell $\Val\setminus\{0\}$. The $\TPas$-morphism
$Z\to\Ord:x\mapsto \ord\, x^2$ is $\cT$-prepared since with the
$\cT$-presentation
$$
\lambda:Z\to \Ord\times \Res_1\times\Val\setminus\{0\}:x\mapsto
(\ord\, x,\ac\, x, x)
$$
and the map $g:\Ord\times \Res_1\to\Ord:(z,\eta)\mapsto 2z$, one has
$f= g \circ \pi \circ \lambda$.
\end{ex}

The following two theorems lay the technical foundations for
analytic motivic integration, the first one to calculate the
integrals, the second one to prove a change of variables formula,
cf.~\cite{CLoes} for the algebraic setting.
\begin{thm}[Cell decomposition]\label{thm:cd} Let $\cT$ be $\TPas$,
$\Tan$, or $\Tans$, and let $\cL_{\cT}$ be respectively $\LPas$,
$\Lan$, or $\Lans$. Let $X$ be a $\cT$-subassignment of $S\times
\Val$ and let $f:X\to R$ be a $\cT$-morphism with $R$ a Cartesian
product of auxiliary sorts, $S$ a $\cT$-assignment. Then there
exists a finite partition of $X$ into $\cT$-cells $Z$ such that each
of the restrictions $f|_{Z}$ is $\cT$-prepared. Moreover, this can
be done in such a way that for each occurring cell $Z$ one can
choose a presentation $\lambda:Z\to Z_{C}$  onto a cell $Z_C$ with
center $c$, such that $c$ is given by an $\cL_\cT^*$-term, where
$\cL_\cT^*$ is defined by (\ref{Lstar}).
\end{thm}

The following is a fundamental result on the term-structure of
definable functions. The statement of Theorem \ref{normal:an} for
fields of the form $k((t))$, uniform in the field $k$ of
characteristic zero, and ideals $I_2=I_1$ was announced in
\cite{CLoes} and will be proved completely here.
\begin{thm}[Term structure of definable morphisms]\label{normal:an} Let $\cT$ be $\TPas$,
$\Tan$, or $\Tans$, and let $\cL_{\cT}$ be respectively $\LPas$,
$\Lan$, or $\Lans$. Let $f:X\to Y$ be a $\cT$-morphism. Then there
exist a $\cT$-parameterization $g:X\mapsto X'$ and a tuple $h$ of
$\cL_\cT^*$-terms in variables running over $X'$ and taking values
in $Y$ such that $f= h\circ g$. (See~(\ref{Lstar}) and
\ref{def:parameterization} for the definitions.)
\end{thm}

The following notion is only needed for the proof of quantifier
elimination in the language $\TPas$, cf.~similar proofs in
\cite{Cohen}, \cite{Denef2}, \cite{pas}, and \cite{pas2}.
\begin{defn}\label{def:strong}
An $\LPas$-definable function $h:X\to \Val$, with $X$ a Cartesian
product of sorts, is called strongly definable if for each
$\Val$-quantifier free $\LPas$-formula $\varphi(v,y)$, with $y$ a
tuple of variables running over arbitrary sorts, $v$ a
$\Val$-variable, and $x$ running over $X$, there
 exists a $\Val$-quantifier free
$\LPas$-formula $\psi(x,y)$ such that
$$
\varphi(h(x),y)
$$
is $\TPas$-equivalent with
$$
\psi(x,y).
$$
\end{defn}

The next Lemma yields $\LPas^*$-terms picking a specific root when
Hensel's Lemma implies that there exists a unique such one.

\begin{lem}\label{lemhens} Let $Z$ be a $1$-$\TPas$-cell with
$\TPas$-presentation ${\rm id}:Z\to Z=Z_{C, \alpha, \xi, c}$ with
$\xi$ taking values in $\Res_{e+1}$ and such that $c=0$ on $C$. Let
$x$ run over $C$, and $y$ over $\Val$. Let $n>0$ and
$f(x,y)=\sum_{i=0}^{n} a_i(x)y^i$ be a polynomial in $y$ with
$\LPas^*$-terms $a_i(x)$ as coefficients, such that $a_n(x)$ is
nowhere zero on $C$. Suppose that for $(x,y)$ in $Z$
$$
\min_i \ord\, a_i(x)y^i= \ord\, a_{i_0}(x)y^{i_0} \ \mbox{ for
some fixed }i_0\geq 1,
$$
and
$$\ord\, f'(x,y)\leq \ord\, t_0^e a_{i_0}(x)y^{i_0-1},$$
and that there exists a $\TPas$-morphism $d:C\to\Val$ whose graph
lies in $Z$ and satisfies
$$
f(x,d(x))=0.
$$
 Then,
$d$ is the unique such morphism and, after a $\TPas$-parametrization
of $C$, $d$ can be given by an $\LPas^*$-term. Moreover, if the
$a_i$ are strongly definable, then the function $d$ is strongly
definable (see Definition \ref{def:strong}).
\end{lem}
\begin{proof}
The uniqueness of $d$ follows from Hensel's Lemma, cf.~\cite{pas}
and \cite{pas2}. Consider the $\TPas$-parametrization
$$g:C\to C\times
\Res_{2e+1}\times \Ord:x\mapsto (x,\ac_{2e+1}\, d(x),\ord\, d(x)).$$
We prove that, piecewise, $d$ can be given by a term after the
$\TPas$-parametrization $g$; at the end we will glue the pieces
together. Note that $\alpha(x)=\ord\, d(x)$. We first prove that
there exists a $\Val$-term $b$ such that $\ord\, b(g(x))=
\alpha(x)$. Let $f_{I_x}(x,y)$ be the polynomial $\sum_{i\in
I_x}a_i(x)y^i$, with
$$
 I_x:=\{i\in\{0,\ldots,n\} :  \ord\, a_i(x)y^i \leq
\ord\, t_0^{2e+1}a_{i_0}(x)y^{i_0}\}.
 $$
 Note that $I_x$ only depends on $x$, since the valuation of $y$ for
 $(x,y)\in Z$
 only depends on $x$.
 We work piecewise to find $b$. First we work on the piece where
$\mathrm{gcd}(i\in I_x)=1$. After partitioning further to ensure
that the quotients
$$
a_i(x)y^i/a_j(x)y^j
$$
have constant order on $Z$ for $i,j\in I_x$,
 one readily verifies that there exists an $\LPas^*$-term $b$ such
that $\ord\, b(g(x))= \alpha(x)$ (for this, the constant symbol
$t_0$ is needed).
 Now work on the part $\mathrm{gcd}(i\in I_x)=\ell$ for some
$\ell>1$. One obtains, by induction on the degree, an
$\LPas^*$-term $h$ such that
 $$\sum_{i\in I_x}a_i(x)h(g(x))^{i/\ell}=0,\quad \ord\,
 h(g(x))=\ell\alpha(x),\quad \mbox{and}\quad \ac_{e+1} h(g(x))=\xi(x)^\ell.$$
 By the conditions of the lemma, $\ell\not=0$ in $\Res_{e+1}$.
Defining the term $b(x,\eta,a)$ as $(h(x,\eta,a))^{1/\ell}$, one
verifies that $\ord\, b(g(x))= \alpha(x)$ for all $x\in C$ with
$\mathrm{gcd}(i\in I_x)=\ell$.

 \par
Now one has  $d(x)=\tau(g(x))$ with $\tau$ the term
$$\tau(x,\eta,a):=b h_{n,e}(\frac{a_0}{a_{i_0}b^{i_0}},\frac{b
a_1}{a_{i_0}b^{i_0}}, \ldots,\frac{b^n a_n}{a_{i_0}b^{i_0}},
\eta\ac_{2e+1}(1/b)).$$

\par
One can glue $s$ pieces together using extra parameters contained in
the definable subassignment $A:=\{\xi\in \Res_1^s :  \sum_i
\xi_i=1\, \wedge\, (\xi_i=0\, \vee\, \xi_i=1)\}$ to index the
pieces, by noting that for each element $a$ in $A$ there exists a
definable morphism $A\to \Val$, given by an $\cL_\cT^*$-term, which
is the characteristic function of $\{a\}$. The fact that $d$ is
strongly definable when the $a_i$ are will be proved in the proof of
Theorem \ref{thm:cd}.
\end{proof}

The following is a refinement of both Theorem 3.1 of \cite{pas}
and Theorem 3.1 of \cite{pas2}, the refinements being the same as
the list of algebraic refinements in the introduction of section
\ref{sec5}.
\begin{thm}\label{thm1:pas} Let $f(x,y)$ be a polynomial in $y$ with
$\LPas^*$-terms in $x=(x_1,\ldots,x_m)$ as coefficients, $x$ running
over a $\TPas$-assignment $S$. Then there exists an integer $\ell$
and a finite partition of $S\times \Val$ into $\TPas$-cells $Z$ with
presentation $\lambda:Z\to Z_C$ such that $Z_C$ has an
$\LPas^*$-term $c$ as center, and such that, if we write
   $$
   f(x,y) = \sum_i a_i(z) (y-c(z))^i,
   $$
for $(x,y)\in Z$ and $(z,y)=\lambda(x,y)$, then we have
   $$
\ord_\infty\, f(x,y) \leq  \min_i \ord_\infty\, t_0^\ell a_i(z)
(y-c(z))^i
   $$
for all $(x,y)\in Z$. Here, $\ord_\infty$ is the extension of $\ord$
by $\ord_\infty(0)=+\infty$. If one restricts to the theory
$\TPas\cup$(VII)$_0$, one can take $\ell=0$ and one can choose cells
whose angular components take values in $\Res_1$.
\end{thm}
\begin{proof}

\par
The equicharacteristic $0$ result induces the analogous result for
big enough residue field characteristic, leaving only finitely many
residue characteristics and ramification degrees to treat
separately.
\par
We give a proof for $\cT=\TPas\cup$(VII)$_p$ with $p>0$. For
$\cT=\TPas\cup$(VII)$_0$ one can use the same proof with $e=0$,
$\ell_0=0$, and $k=1$.
\par
Let $f$ be of degree $d$ in $y$ and proceed by induction on $d$.
Let $f'(x,y)$ be the derivative of $f$ with respect to $y$.
Applying the induction hypothesis to $f'$, we find a partitioning
of $S\times\Val$ into cells. By replacing $S$ we may suppose that
these cells have the identity mapping as presentation.

\par
First consider a $0$-cell $Z=Z_{C, c}$ in the partition. Then we can
write $ f(x,y)=f(x,c(x))=\tau(x)$ for some $\LPas^*$-term $\tau$,
for $(x,y)\in Z$, and the theorem follows.

\par
Next consider a $1$-cell $Z=Z_{C, \alpha, \xi, c}$ in this
partition. Write $f(x,y) = \sum_i a_i(x) (y-c(x))^i$ for $(x,y)\in
Z$. There is some $\ell_0$ such that for all $(x,y)\in Z$
   $$
\ord\, f'(x,y) \leq \min_i \ord\, t_0^{\ell_0}ia_i(x)
(y-c(x))^{i-1}.
   $$
We may suppose that $a_i(x)$ is either identically zero or else
never zero on $C$ for each $i$. Put $I=\{i :  \forall x\in C\
a_i(x)\not=0\}$ and $J=\{(i,j)\in I\times I :  i>j\}$. We may
suppose that $J$ is nonempty since the case $J=\emptyset$ is
trivial. Put $\Theta:=\{<,>,=\}^J$. For
$\theta=(\sq_{ij})\in\Theta$, put
 $$C_\theta= \{x\in C :  \forall (i,j)\in J\ i\alpha(x)+\ord\, a_i(x)
 \sq_{ij}j\alpha(x)+\ord\, a_j(x)\}.
 $$
Ignoring the $C_\theta$ which are empty, this gives a partition of
$C$ and hence a partition of $Z$ into cells
$Z_\theta=Z_{C_\theta,\alpha|_{C_\theta}, \xi|_{C_\theta},
c|_{C_\theta}}$.
 Fix $\theta\in \Theta$.  We may suppose that
$Z=Z_\theta$.
 The case that $\ord\, a_0(x)<i\alpha(x)+ \ord\, a_i(x)$ for all
$i\geq 1$ and all $x\in C$ follows trivially. Hence, we may suppose
that there exists $i_0>0$ such that $i_0\alpha(x) +\ord\,
a_{i_0}(x)\leq i\alpha(x)+ \ord\, a_i(x)$ for all $x\in C$ and all
$i$.
 Put $e:=\ell_0+\ord\, i_0$. Let $\xi$ take values in
$\Res_k$. By either enlarging $e$ and $\ell_0$ or enlarging $k$, we
may suppose that $k=e+1$.
 Define the subassignments $B_i\subset Z$ by
\begin{gather*}
B_{1}:=  \{(x,y)\in Z :  \\
 \forall\, z\,  \big(
\ord(z-y)>\alpha(x)+\ord\, t_0^e\, \to\, \ord f(x,z) \leq
i_0\alpha(x) +\ord\, t_0^{2e} a_{i_0}(x)\big)\}
 \end{gather*}
 and $B_{2}:= Z\setminus B_1.$
 If $B_{i}$ is nonempty, it is equal to the cell $Z_{C_{
i},\alpha|_{C_{i}}, \xi|_{C_{i}}, c|_{C_{i}}}$ for some
$\TPas$-definable assignment $C_{i}\subset C$. Moreover, the $B_i$
can be described without using new $\Val$-quantifiers, by using the
maps $\ac_m$ for big enough $m$. On $B_{1}$ the theorem holds with
$\ell=2e$. On $B_2$, by Lemma \ref{lemhens}, there exists a unique
definable function $d:C_2\to \Val$ such that, for each $x$ in $C_2$,
$(x,d(x))$ lies in $B_{2}$ and $f(x,d(x))=0$. Again by Lemma
\ref{lemhens}, we may suppose that $d(x)$ is given by an
$\LPas^*$-term.
 Let $D$ be the $\TPas$-definable assignment $\{j\in\Ord :  j>0\}$.
 For
\begin{gather*}
C_2':=C_2\times D\times (\Res_1\setminus\{0\}),\\
 \beta:C_2'\to
\Ord: (x,j,z)\mapsto \alpha(x)+e+j,\\
\eta:C_2'\to \Res_1\setminus\{0\}: (x,j,z)\mapsto z,\\
d':C_2'\to \Val:(x,j,z)\mapsto d(x),\\
\end{gather*}
consider the cell $Z_{C_2', \beta, \eta, d'}$ and its projection
$\pi$ to $C_2\times \Val$. One checks that $B_{2}$ is the disjoint
union of the $1$-cell $\pi(Z_{C_2', \beta, \eta, d'})$, with
presentation $\pi^{-1}$, and the $0$-cell $Z_{C_2,d}$. Moreover, if
one writes $f(x,y)= \sum b_i(z)(y-d'(z))^i$ for $(z,y)$ in $Z_{C_2',
\beta, \eta, d'}$ and $(x,y)=\pi(z,y)$, one has $ \ord\, f(x,y) =
\ord\, b_1(z) (y-d'(z))$, which can be seen using a Taylor expansion
of $f$ around $d'$.
 This finishes the proof.
\end{proof}

\begin{rem}\label{rem:ideals}
\item[(i)] By Theorem \ref{thm1:pas}, one can probably add arbitrary
angular components modulo $I_j$ for a collection of nonzero ideals
${I_j}_{j\in J}$ to $\LPas$, $\Tan$, or $\Tans$. When one enlarges
the language $\LL_{\rm Res}$ to the the full induced language,
(which can be richer in the analytic than in the algebraic case),
one can probably obtain a form of quantifier elimination and cell
decomposition. Most likely, one also gets a similar term structure
result (even without introducing new $(\cdot,\cdot,\cdot)_j^{1/m}$
or $h_{m,j}$ for the new ideals). This should follow from Theorem
\ref{thm1:pas}.
\item[(ii)] Similar proofs should hold to show that, if one restricts to the theory
$\cT\cup$(VII)$_0$ in Theorem \ref{normal:an}, one can take for $h$
in Theorem \ref{normal:an} a tuple of $\cL_\cT^\diamond$-terms, with
$\cL_\cT^\diamond$ the language
 \begin{equation*}
\cL^\diamond:=\cL\cup_{m>0}\{(\cdot,\cdot,\cdot)_0^{1/m},\
h_{m,0}\}.
 \end{equation*}
\end{rem}

\subsection*{Proof of Theorems \ref{thm:cd} and \ref{thm:qe:an} for $\TPas$}
First suppose that $X=X_0=\Val^{m+1}$ and that $f=f_0$ is the map
\begin{equation}\label{cd:polynomap}
f_0:\Val^{m+1}\to \Res_n^\ell\times \Ord^\ell: x\mapsto
(\ac_n(g_i(x,t)),\ord(g_i(x,t)))_i,
 \end{equation}
with $g_i(x_1,\ldots,x_m,t)$ polynomials over $\ZZ$, $m\geq 0$,
$\ell,n>0$, $i=1,\ldots,\ell$.
 By Theorem \ref{thm1:pas}, the result
for $\ell=1$ follows rather immediately. It is from this partial
result for $\ell=1$ that one deduces the final statement of Lemma
\ref{lemhens} in the same way as this is proved in \cite{pas} and
\cite{pas2}. We will not recall this proof of the final statement of
Lemma \ref{lemhens}.
\par
By induction on $\ell$, we may suppose that the result holds for
$G_1:=(\ac_n\, g_i,\ord \,g_i)_{i=1}^{\ell-1}$ and for
$G_2:=(\ac_n\,g_\ell,\ord\,g_\ell)$. This gives us two finite
partitions $\{Z_{ij}\}$ such that $G_i$ is prepared on $Z_{ij}$ for
each $j$ and $i=1,2$.
 Choose $Z_1:=Z_{1j}$ and $Z_2:=Z_{2j'}$. It is enough to partition
$Z_1\cap Z_2$ into cells such that $f_0$ is prepared on these cells.
If $Z_1$ or $Z_2$ is a $0$-cell, this is easy, so we may suppose
that $Z_1$ is a $1$-cell with presentation $\lambda_1:Z_1\to
Z_{C_1}=Z_{C_1, \alpha_1, \xi_1, c_1}$ and $Z_2$ a $1$-cell with
presentation $\lambda_2:Z_2\to Z_{C_2}=Z_{C_2, \alpha_2, \xi_2,
c_2}$. We may suppose that $\pi(Z_1)=\pi(Z_2)$ with $\pi:X\to S$ the
projection, that $\xi_1$ and $\xi_2$ take values in $\Res_k$, and
that $Z_{C_i}\subset Z\times R$ with $R$ a fixed product of
auxiliary sorts for $i=1,2$. Under these suppositions it follows
from the non archimedean property that $Z_1\cap Z_2$ is already a
cell on which the function $f_0$ is $\TPas$-prepared, where one can
use the presentation
  \begin{gather*}\lambda_{12}:Z_1\cap Z_2\to \lambda_{12}(Z_1\cap Z_2)\subset Z\times
 R\times \{0,1\}:\\ (x,t)\mapsto \left\{\begin{array}{ll}(\lambda_{1}(x,t),0) & \mbox{if }
  \alpha_1\geq  \alpha_2,\\
(\lambda_{2}(x,t),1) & \mbox{else,}
\end{array}\right.
\end{gather*}
and where we write $\alpha_i$ for $\alpha_i(\lambda_i(x,t))$. Here,
$\lambda_{12}(Z_1\cap Z_2)$ has as center $c_1d_0+c_2d_1$, where
$d_i$ is the $\LPas^*$-term from $\Res_1$ to $\Val$ which is the
characteristic function of $\{i\}$ for $i=0,1$.

\par
By exploiting the proof of this partial result for general $\ell$,
one can ensure that all occurring centers are strongly definable and
that there are no $\Val$-quantifiers introduced in the process of
the cell decomposition. From this partial result for general $\ell$
one deduces quantifier elimination for $\TPas$ in the language
$\LPas$ as in \cite{pas} and \cite{pas2}.

\par
Now let $f:X\to R$ be a general $\cT$-morphism with $R$ a Cartesian
product of auxiliary sorts and $X$ an arbitrary $\cT$-assignment.
Let $f_1,\ldots,f_t$ be all the polynomials in the ${\rm
Val}$-variables, say, $x_1,\ldots,x_{m+1}$ occurring in the formulas
describing the $X$ and $f$, where we may suppose that these formulas
do not contain quantifiers over the valued field sort. Apply the
above case of cell decomposition to the polynomials $f_i$. This
yields a partition of $\Val^{m+1}$ into cells $Z_i$ with
presentations $\lambda_i:Z_i\to Z_{C_i}$ and centers $c_i$. Write
$x=(x_1,\ldots,x_{m+1})$ for the ${\rm Val}$-variables,
$\xi=(\xi_j)$ for the ${\rm Res}$-variables and $z=(z_j)$ for the
${\rm Ord}$-variables on $Z_{C_i}$.
 If $Z_i$ is a $1$-cell, we may suppose that for $(x,\xi,z)$ in
$Z_{C_i}$ we have $\ord(x_{m+1}-c_i)=z_1$ and
$\ac_n(x_{m+1}-c_i)=\xi_1$, by changing the presentation of $Z_i$
if necessary (that is, by adding more ${\rm Ord}$-variables and
${\rm Res}$-variables). By changing the presentation as before if
necessary, we may also assume that
\[
\begin{array}{l}
\ord f_j(x)= z_{k_j},\\
\ac_n f_j(x) =  \xi_{l_j},
\end{array}\]
for $(x,\xi,z)$ in the $1$-cell $Z_{C_i}$, where the indices $k_j$
and $l_j$ only depend on $j$ and $i$.

\par
 Since the condition $f(x)=0$ is equivalent to $\ac_n(f(x))=0$, we
may suppose that, in the formulas describing $X$ and $f$ the only
terms involving ${\rm Val}$-variables are of the forms $\ord\,
f_j(x)$ and $\ac_n\, f_j(x)$. Combining this with the above
description of $\ord f_j(x)$ and $\ac_n f_j(x)$, one sees that the
value of $f$ only depends on variables running over the bases of the
cells. Hence, $f$ is $\TPas$-prepared on these cells. 
\begin{flushright}
$\square$
\end{flushright}

\subsection*{Proof of Theorem \ref{thm:cd} for $\Tan$ and $\Tans$}
Let $\cM$ be a model of $\cT$, $a$ a $\Val$-tuple of $\cM$, $\cM_a$
the $\cL_\cT$-substructure of $\cM$ generated by $a$, $K_a$ the
valued field of $\cM_a$, $\cL_\cT(K_a)$ the language $\cL_\cT$
together with constant symbols for the elements of $K_a$, and
$\cT(K_a)$ the $\cL_\cT(K_a)$-theory of all models of $\cT$
containing the structure $\cM_a$.  Let $\LPas(K_a)$ be the language
$\LPas$ together with constant symbols for the elements of $K_a$ and
$\TPas(K_a)$ the $\LPas(K_a)$-theory of all $\TPas$-models
containing the structure $\cM_a$.

\par
First we prove some special cases of Theorem \ref{thm:cd} for the
theory $\cT(K_a)$. Suppose first that $X=\Val$ and that $f$ is the
map
$$
f:X\to \Res_n^\ell\times \Ord^\ell: x\mapsto
(\ac_n(g_i(x)),\ord(g_i(x)))_i,
$$
with the $g_i$ $\cL_\cT(K_a)$-terms in the variable $x$ for
$i=1,\ldots,\ell$. In the case that $\cT$ is $\Tan$, apply Theorem
\ref{lem:toholom} to the terms $g_i$ and to the terms $g_i(x^{-1})$.
In the case that $\cT$ is $\Tans$, there is a $\cL_\cT$-term which
presents a valued field element with minimal positive valuation by
(IV)$_T$, hence we can apply Corollary \ref{cor:strictrational} to
the terms $g_i$ and to the terms $g_i(x^{-1})$.
 In this way we find a finite partition of
$X_0:=\{x\in\Val :  \ord(x)\geq 0\}$ into
$\LPas(K_a)$-$\cT(K_a)$-assignments $X_j$ given by $K_a$-annulus
formulas $\varphi_j$, rational functions $h_{ij}(x)$ with
coefficients in $K_a$, a polynomial $F(x)$ over $K_a$, and very
strong units $U_{ij}\in \cO_{K_a}(\phi_j)$, such that for all $i,j$
and all $x\in X_j$
$$
F(x)=0\ \vee\ g_i(x) = U_{ij}(x) h_{ij}(x),
$$
where we mean by an $\LPas(K_a)$-$\cT(K_a)$-assignment a
$\cT(K_a)$-assignment which can be defined by an
$\LPas(K_a)$-formula. If $I_1\not= I_2$, the separated analytic
structure collapses to a strictly convergent analytic structure, and
thus, by Corollary \ref{cor:strictrational}, we can even assume that
$\ac_n(U_{ij})=1$ on $X_j$. Up to the transformation $x\mapsto
x^{-1}$, we can partition $X_1:=\{x\in\Val :  \ord(x)< 0\}$ in a
similar way.
 Now apply Theorem \ref{thm:cd} for the theory $\TPas(K_a)$ to the
$\LPas(K_a)$-$\cT(K_a)$-assignments $X_{jk}$ and the functions $x\in
X_{jk}\mapsto (\ac(F(x)),\ac_n\, h_{ij}(x),\ord\, h_{ij}(x))_{i}$ to
refine the partition and to finish the proof of $X$ and $f$ of the
above form.

 \par
Next we suppose that $X=\Val$ and $f$ is an arbitrary
$\cT(K_a)$-morphism $f:X\to \Res_n^\ell\times \Ord^\ell$. Apply
Theorem \ref{thm:qe:an} to obtain a formula $\varphi$ without
$\Val$-quantifiers, as in (\ref{form}), which describes the graph of
$f$. Then, let $g_1,\ldots,g_r$ be the $\cL_\cT(K_a)$-$\Val$-terms
occurring in $\varphi$. Applying the previous case to the terms
$g_i$, the case of this $f$ easily follows, cf. the analogous step
in the proof of Theorem \ref{thm:cd} for $\TPas$.

\par
Next we suppose that $X=\Val^{m+1}$ and that $f$ is an arbitrary
$\cT$-morphism $f:X\to \Res_n^\ell\times \Ord^\ell$, $\ell,m>0$. In
this case, the theorem is reduced by a compactness argument to the
case $m=0$, as follows. Suppose that for every candidate $\cT$-cell
decomposition of $X$ into $\cT$-cells $A_i$, with $\cL_\cT^*$-terms
as centers, and $\cT$-prepared functions $g_i:A_i\to
\Res_n^\ell\times \Ord^\ell$, this data is not the data of a cell
decomposition of $X$ which prepares $f$. This is equivalent to
saying that for each such candidate cell decomposition there exists
a model with valued field $K$ (with $A$-analytic structure) and
$a\in K^m$ such that either the fibers of the $\cT$-cells $A_i$
above $a$ (under the projection $\Val^{m+1}\to \Val^m$) are not a
$\cT(a)$-cell decomposition of $\Val$, or the fibers of the
functions $g_i$ above $a$ do not coincide with $f$ on the fiber of
$A_i$ above $a$. Then, by compactness, there exists a model with
valued field $K'$ and $a'\in K'{}^m$ such that $\Val$ can not be
partitioned into $\cT(K_a)$-cells on which the fibers of the
functions $g_i$ above $a$ are prepared, which contradicts the
previous case for $X=\Val$. Moreover, this construction ensures that
we can work with $\cL_\cT^*$-terms as centers of the cells.

\par
Finally, the general Theorem follows from this case similarly as the
general case is obtained in the proof of Theorem \ref{thm:cd} for
$\TPas$.
\begin{flushright} $\square$
\end{flushright}

\subsection*{Proof of Theorem \ref{normal:an}}
By working componentwise it is enough to prove the theorem for
$Y=\Val$. Let ${\rm Graph}(f)\subset X\times \Val$ be the
$\cT$-assignment which is the graph of $f$, and suppose it is
described by an $\cL_\cT$-formula $\varphi$ in the form given by the
quantifier elimination Theorem \ref{thm:qe:an}; let $g_j$ be the
$\cL_\cT$-terms occurring in this formula. Apply Theorem
\ref{thm:cd} for $\cT$ to the terms $g_i$. Doing so, all occurring
centers of the cells are given by $\cL_\cT^*$-terms. For each
occurring cell $Z_i$, let $Z'_i$ be $\lambda_i^{-1}({\rm
Graph}(c_i))\cap {\rm Graph}(f)$, where $\lambda_i$ is the
presentation of $Z_i$ and $c_i$ its center. Clearly each $Z_i'$ is a
$0$-cell with presentation the restriction of $\lambda_i$ to $Z_i'$.
It follows from the special form of $\varphi$ (as given by the
application of Theorem \ref{thm:qe:an}) that the cells $Z_i'$ form a
cell decomposition of ${\rm Graph}(f)$ and one concludes that the
restriction of $f$ to each of finitely many pieces in a partition of
$X$ satisfies the statement. Now the Theorem follows by gluing the
pieces together using extra parameters as in the proof of Lemma
\ref{lemhens}.
\begin{flushright}
$\square$
\end{flushright}

\section{Applications to analytic motivic integration}
\label{part:motivic}

Let $\cO_F$ be the ring of integers of a number field $F$. Let
$\cA_F$ be the class of all finite field extensions of all $p$-adic
completions of $F$, and $\cB_F$ the class of all local fields of
positive characteristic which are rings over $\cO_F$. For a fixed
prime $p$ and integer $n>0$, let $\cA_{F,p,n}$ be the subset of
$\cA_F$ consisting of all fields with residue field of
characteristic $p$ and with degree of ramification fixed by
$\ord_p(p)=n$.

\par
For $K\in \cA_F\cup\cB_F$ write $K^\circ$ for the valuation ring,
$\pi_K$ for a uniformizer, $\Kt$ for the residue field, and $q_K$
for $\sharp \Kt$. By $T_m(\cO_F[[t]])$ denote the ring of strictly
convergent power series in $m$ variables over $\cO_F[[t]]$. For each
$K$ in $\cA_F\cup \cB_F$ and each power series
$f=\sum_{i\in\NN^m}a_i(t)X^i$ in $T_m(\cO_F[[t]])$ define the
analytic function
$$
f_K:(K^\circ)^m\to K^\circ:x\mapsto \sum_{i\in\NN^m}a_i(\pi_K)x^i,
$$
and extend this by zero to a function $K^{m}\to K$.

\par
In the terminology of section \ref{part:an}, we have thus fixed the
strictly convergent analytic $\cO_F[[t]]$-structure on all the
fields $K\in\cA_F\cup \cB_F$.

\par
Let $\cL^F$ be the language $\cL_{T(\cO_F[[t]])}$ with $\LL_{\Ord}$
the Presburger language $\LL_{\rm Pres}=(+,-,0,1,\leq,\{\equiv\bmod
n\}_n)$ and $\LL_{\Res}$ the language $\LL_{{\rm Res},0}$
(cf.~section \ref{sec:language}). Define the $\cL^F$-theory $\cT^F$
as $\Tans$ together with the axiom $t_0\not=1$ (that is, with higher
order angular components, see section \ref{sec:language}), and
axioms describing the congruence relations modulo $n$ in the natural
way.

\par
Let $W$ be an $\cL^F$-formula with $m$ free valued field variables
and no other free variables. (Note that $W$ determines a
$\cT^F$-assignment in the sense of section \ref{sec:assignm}, but
this is not needed here.) For each $K\in\cA_F\cup\cB_F$, we obtain a
set $W_K\subset K^m$ by interpreting the formula $W$ in the natural
way. In a similar way, a $\cT^F$-morphism $f$ from $W$ to the valued
field\footnote{By this we mean an $\cL^F$-formula $\varphi$ such
that the set described by $\varphi$ in any model $\cM$ of $\cT^F$ is
the graph of a function from the $\cM$-rational points on $W$ to the
valued field of $\cM$, cf.~section \ref{sec:assignm}.} determines a
function $f_K:W_K\to K$.

\par
Suppose now that the set $W_K$ is contained in $(K^\circ)^m$ for
each $K\in\cA_F\cup\cB_F$. Fix $\cT^F$-morphisms $f_1$ and $f_2$
from $W$ to the valued field, such that the images of the $f_{iK}$
lie in $K^\circ$ for each $K\in\cA_F\cup\cB_F$.

\par
For each $K\in\cA_F\cup\cB_F$ and $s\geq0$ a real number, we
consider
\begin{equation}\label{eq:aks}
a_K(s):= \int_{W_K} |f_{1K}(x)|^s |f_{2K}(x)||dx|.
\end{equation}

\par
It is well known by work of Denef and van den Dries \cite{DvdD}
that, for each fixed $K\in \cA_F$, $a_K(s)$ is a rational function
in $q_K^{-s}$. We prove that also for fixed $K\in \cB_F$ with big
enough characteristic, $a_K(s)$ is a rational function in
$q_K^{-s}$. Moreover, we give a geometric meaning to how the $a_K$
vary for $K\in \cA_F$, and, when the characteristic is big enough,
also for $K\in \cB_F$.

\par
Let ${\rm Var}_{\cO_F}$ denote the collection of isomorphism classes
of algebraic varieties over $\cO_F$ and let ${\rm Form}_{\cO_F}$ be
the collection of equivalence classes of formulas\footnote{Two
formulas are equivalent in this language if they have the same
$R$-rational points for every ring $R$ over $\cO_F$.} in the
language of rings with coefficients in $\cO_F$. Define the rings
   $$
\cM(\Var_{\cO_F}):= \QQ[T,T^{-1},{\rm
Var}_{\cO_F},\frac{1}{\AA^1_{\cO_F}},\{\frac{1}{1-\AA^{b}_{\cO_F}T^a}
\}_{(a,b) \in J} ]
   $$
and
   $$
\cM({\rm Form}_{\cO_F}):= \QQ[T,T^{-1},{\rm
Var}_{\cO_F},\frac{1}{\AA^1_{\cO_F}},\{\frac{1}{1-\AA^{b}_{\cO_F}T^a}
\}_{(a,b) \in J} ],
   $$
with $J=\{(a,b)\in\ZZ^2 :  a\geq 0,\ b<0\}$, and where we write
$\AA^\ell_{\cO_F}$ for the isomorphism class of the formula
$x_1=x_1\wedge\ldots\wedge x_\ell=x_\ell$ (which has the set
$R^\ell$ as $R$-rational points for any ring $R$ over $\cO_F$),
$\ell\geq 0$.

\par
For each finite field $k$ over $\cO_F$ with $q_k$ elements, we write
${\rm Count}_{k}$ for the ring morphisms
   $$
{\rm Count}_{k}:\cM(\Var_{\cO_F}) \to
\QQ[q_k^{-s},q_k^s,\{\frac{1}{1-q_k^{-as+b}} \}_{(a,b) \in J} ]
   $$
   and
 $$
{\rm Count}_{k}:\cM({\rm Form}_{\cO_F}) \to
\QQ[q_k^{-s},q_k^s,\{\frac{1}{1-q_k^{-as+b}} \}_{(a,b) \in J} ]
   $$
which send $T$ to $q_k^{-s}$, $Y\in{\rm Var}_{\cO_F}$ to the
number of $k$-rational points on $Y$ and $\varphi \in {\rm
Form}_{\cO_F}$ to the number of $k$-rational points on $\varphi$.

\par
We prove the following generalization of Theorem \ref{thm:int}:
\begin{thm}\label{thm:int:mets}
\item[(i)] There exists a (non-unique) element $X\in
\cM(\Var_{\cO_F})$ and a number $N$ such that for each
$K\in\cA_F\cup \cB_F$ with ${\rm Char} \widetilde K>N$ one has
$$
a_K(s)= {\rm Count}_{\widetilde K}(X).
$$
In particular, for $K\in\cA_F\cup \cB_F$ with ${\rm Char}
\widetilde K$ big enough, $a_K(s)$ only depends on $\widetilde K$.

 \item[(ii)] For fixed prime $p$ and $n>0$
there exists a (non-unique) element $X_{p,n}\in \cM({\rm
Form}_{\cO_F})$ such that for each $K\in\cA_{F,p,n}$ one has
$$
a_K(s)= {\rm Count}_{\widetilde K}(X_{p,n}).
$$
\end{thm}

\begin{proof}[Proof of Theorem \ref{thm:int:mets}.]
The Cell Decomposition Theorem \ref{thm:cd} together with the
Quantifier Elimination Theorem \ref{thm:qe:an} translates the
calculation of the $a_K(s)$ in a nowadays standard way into
calculations of the form
$$
\sum_i {\rm Count}_{\widetilde K}(\varphi_i)\sum_{z\in S_i\subset
\ZZ^m}q_K^{-\alpha_i(z) - s\beta_i(z)}
$$
with $\alpha_i$, $\beta_i:S_i\to\NN$ Presburger functions on the
Presburger sets $S_i$, and the $\varphi_i$ $\LL_{\rm
ring}$-formulas\footnote{Here, we use that, for any $\LL_{{\rm
Res}}$-formula $\varphi$, there exists an $\LL_{\rm ring}$-formula
$\psi$, such that for all $K\in \cA_F\cup \cB_F$ the number of
$K^\circ/(\pi_K^m\Ko)$-rational points on $\varphi$ is the same as
the number of $\Kt$-rational points on $\psi$. For this, use the
definable bijection $K^\circ/(\pi_K^m\Ko)\to\Kt^m$.}. One such
expression works for $K\in\cA_F\cup\cB_F$ with residue field
characteristic big enough, and one needs another such expression for
each class $\cA_{n,p}$. Using Lemma~3.2 of \cite{D85} on the
summation of such Presburger functions and their exponentials, one
obtains a number $N_0$ and objects $X_0$, $X_{p,n}$ of $\cM({\rm
Form}_{\cO_F})$ such that
\begin{itemize}
 \item[(i)] for each $K\in\cA_F\cup\cB_F$
with ${\rm Char} \widetilde K>N_0$ one has
$$
a_K(s)= {\rm Count}_{\widetilde K}(X),
$$
 \item[(ii)] for each $K\in\cA_{F,p,n}$ one has
$$
a_K(s)= {\rm Count}_{\widetilde K}(X_{p,n}).
$$
\end{itemize}
Associate virtual motives $Y_i$ to the (isomorphism classes of)
formulas $\varphi_i$ occurring in $X_0$ using Theorem 6.4.1 of
\cite{DL3}. Since such a virtual motive lives in a ring of virtual
motives which is the image of (a certain localisation of) the
Grothendieck ring of varieties over $F$ (cf.~Theorem 6.4.1 of
\cite{DL3}), one can find $A_i\in \cM(\Var_{\cO_F})$ such that the
number of $k$-points on $Y_i$ is equal to ${\rm Count}_{k}(A_i)$ for
all finite fields $k$ of characteristic $>N_1$ for some $N_1$. One
then easily finds $X\in \cM(\Var_{\cO_F})$ for which the Theorem
holds for $N=\max(N_0,N_1)$.
\end{proof}

\begin{rem}
The extension of the Denef-Pas cell decomposition to the theory
$\TPas$ which allows for $I_2\not =I_1$, i.e.~, which allows for
higher order angular components $\ac_n$ also in equicharacteristic
zero, gives a new view on the foundational work by Denef and Loeser
on geometric motivic integration \cite{DL2} and subsequent work.
More precisely, using $\ac_n$ to define a broader but analogous
class of semialgebraic sets than in \cite{DL2}, stable sets,
cylinders, and weakly stable sets would be semialgebraic, which is
not the case in \cite{DL2}.
\end{rem}


\begin{thebibliography}{10}

\bibitem{bgr}S.~Bosch, U.~G\"untzer and R.~Remmert \emph{Non-archimedean
Analysis}, Springer-Verlag (1984).

\bibitem{Celikler} Y.~F.~\c{C}elikler \emph{Dimension theory and
parameterized normalization for D-semianalytic sets over
non-Archimedean fields}, to appear in J. Symbolic Logic.

\bibitem{c} R.~Cluckers, \emph{Analytic $p$-adic cell decomposition
and integrals}, Trans. Amer. Math. Soc., {\bf 356} (2004) 1489 -
1499, math.NT/0206161.


\bibitem{CLoes} R.~Cluckers, F.~Loeser \emph{Constructible motivic
functions and motivic integration}, math.AG/0410203.

\bibitem{CLoes1} R.~Cluckers, F.~Loeser \emph{Fonctions constructible
et int\'egration motivic I}, Comptes rendus de l'Acad\'emie des
Sciences, {\bf 339} (2004)  411 - 416 math.AG/0403349.

\bibitem{CLoes2} R.~Cluckers, F.~Loeser \emph{Fonctions constructible
et int\'egration motivic II}, Comptes rendus de l'Acad\'emie des
Sciences, 339 (2004)  487 - 492, math.AG/0403350.

\bibitem{Cohen} P.~J.~Cohen \emph{Decision procedures for real and
$p$-adic fields}, Comm. Pure Appl. Math., \textbf{22} (1969)
131-151.

\bibitem{Denef}
J.~Denef \emph{The rationality of the {P}oincar\'e series
associated to the $p$-adic points on a variety}, Inventiones
Mathematicae {\bf 77} (1984) 1-23.


\bibitem{D85}
J. Denef, \textit{On the evaluation of certain $p$-adic
integrals}, S\'eminaire de th\'eorie des nombres, Paris 1983--84,
Progr. Math., Birkh\"auser Boston, Boston, MA \textbf{59} (1985)
25--47.


\bibitem{Denef2}
J.~Denef \emph{$p$-adic semialgebraic sets and cell
decomposition}, Journal f{\"u}r die reine und angewandte
Mathematik {\bf 369} (1986) 154-166.


\bibitem{DL2} J.~Denef, F.~Loeser,
\emph{Germs of arcs on singular algebraic varieties and motivic
integration}, Inventiones Mathematicae, {\bf 135} (1999), 201-232.

\bibitem{DL} J.~Denef, F.~Loeser,
\emph{Definable sets, motives and $p$-adic integrals}, Journal of
the  American Mathematical Society, {\bf 14} (2001) no.~2,
429-469.

\bibitem{DL3} J.~Denef, F.~Loeser, \emph{On some rational generating series occurring in
arithmetic geometry},  In Geometric Aspects of Dwork Theory,
edited by A. Adolphson, F. Baldassarri, P. Berthelot, N. Katz and
F. Loeser, \textbf{1}, de Gruyter (2004)  509-526,
math.NT/0212202.


\bibitem{DvdD} J.~Denef, {L.~van den}~Dries, \emph{$p$-adic and real
subanalytic sets}, Annals of Mathematics, {\bf 128} (1988)
79--138,

\bibitem{vdd1}L.~van den Dries \emph{Analytic Ax-Kochen-Ershov
theorems} Contemporary Mathematics, {\bf 131} (1992) 379-398.

\bibitem{vdd2} L.~van~den~Dries, notes on cell decomposition.

\bibitem{dhm} L.~van~den~Dries, D.~Haskell and D.~Macpherson,
\emph{One-dimensional $p$-adic subanalytic sets}, J. London Math.
Soc. (2), {\bf 59} (1999) 1-20.

\bibitem{endler} O.~Endler, \emph{Valuation Theory}, Springer-Verlag, 1972.

\bibitem{fp} J.~Fresnel and M.~van~der~Put, \emph{G\'eom\'etrie Analytique
Rigide et Applications}, Birkh\"auser (1981).

\bibitem{FP} J.~Fresnel and M.~van~der~Put, \emph{Rigid Geometry and Applications,} Birkh\"auser (2004).

\bibitem{Kazhdan} D.~Kazhdan, \emph{An algebraic integration}  Mathematics:
frontiers and perspectives, Amer. Math. Soc., Providence, RI,
(2000)  93--115.

\bibitem{Kuhl} F.-V.~Kuhlmann, \emph{Quantifier elimination for {H}enselian fields relative to additive and multiplicative
congruences} Israel Journal of Mathematics \textbf{85} (1994)
277-306.

\bibitem{lang} S.~Lang, \emph{Algebra}, Addison-Wesley.

\bibitem{LR} J.-M.~Lion,  J.-P.~Rolin, \emph{Int{\'e}gration des fonctions sous-analytiques
et volumes des sous-ensembles sous-analytiques} (French)
[Integration of subanalytic functions and volumes of subanalytic
subspaces], Ann. {I}nst. {F}ourier, \textbf{48} (1998) no. 3,
755--767.

\bibitem{lp} L.~Lipshitz and Z.~Robinson, \emph{Rigid subanalytic subsets
of the line and the plane}, Amer. J. Math., {\bf 118}(1996)
493-527.

\bibitem{cs} L.~Lipshitz and Z.~Robinson, \emph{Rigid subanalytic subsets
of curves and surfaces}, J. London Math. Soc. (2), {\bf 59} (1999)
895-921.

\bibitem{LRSep} L.~Lipshitz and Z.~Robinson, \emph{Rings of separated power series}, Ast\'erisque \text{264} (2000) 3-108.

\bibitem{mod} L.~Lipshitz and Z.~Robinson, \emph{Model completeness and subanalytic sets},
Ast\'erisque, \text{264} (2000) 109-126.

\bibitem{unifrid} L.~Lipshitz and Z.~Robinson, \emph{Uniform properties of
rigid subanalytic sets}, to appear in Trans. Amer. Math. Soc.,
available at www.math.purdue.edu/$\sim$lipshitz/.

\bibitem{SB1} F.~Loeser, J.~Sebag, \emph{Motivic integration on smooth rigid
varieties and invariants of degenerations} Duke Math. J.
\textbf{119} (2003) no. 2, 315--344.

\bibitem{MacintRat} A.~Macintyre,  \emph{Rationality of $p$-adic Poincar\'e series:
uniformity in $p$}, Ann. Pure Appl. Logic \textbf{49} (1990) no.
1, 31-74.

\bibitem{NS} J.~Nicaise, J.~Sebag, \emph{Invariant de Serre et fibre de Milnor
analytique}, (French) [The Serre invariant and the analytic Milnor
fiber] available at
www.wis.kuleuven.ac.be/algebra/artikels/artikelse.htm.

\bibitem{pas} J.~Pas, \emph{Uniform $p$-adic cell decomposition and local
zeta-functions}, J. Reine Angew. Math., {\bf 399}(1989) 137-172.

\bibitem{pas2} J.~Pas, \emph{Cell decomposition and local
zeta-functions in a tower of unramified extensions of a $p$-adic
field}, Proc. London Math. Soc., {\bf 60}(1990) 37-67.

\bibitem{Scanlon} T.~Scanlon, \emph{Quantifier elimination for the relative Frobenius} in Valuation
Theory and Its Applications Volume II, conference proceedings of
the International Conference on Valuation Theory (Saskatoon,
1999), Franz-Viktor Kuhlmann, Salma Kuhlmann, and Murray Marshall,
eds., Fields Institute Communications Series, (AMS, Providence),
2003, 323 - 352.


\bibitem{Seb2}
J. Sebag \textit{Rationalit\'e des s\'eries de Poincar\'e et des
fonctions z\^eta motiviques} (French) [Rationality of Poincare
series and motivic zeta functions]  Manuscripta Math. \textbf{115}
(2004), no. 2, 125--162.


\bibitem{SB2} J.~Sebag, \emph{Int\'egration motivique sur les sch\'emas
formels} (French) [Motivic integration on formal schemes]
Bull.~Soc.~Math.~France  \textbf{132}  (2004) no. 1, 1--54.



\end{thebibliography}
\end{document}